\documentclass[preprint,11pt]{elsarticle}

\usepackage{lineno,hyperref}
\usepackage{graphicx}
\usepackage{{ntheorem}}
\usepackage{amssymb}
\usepackage{multicol}
\usepackage{anyfontsize}
\usepackage[T1]{fontenc}
\usepackage{calligra}
\usepackage[x11names,table]{xcolor}
\usepackage{slashbox}
\usepackage{multirow}
\usepackage{pdfpages}
\usepackage{nicefrac}
\usepackage{amsmath}
\usepackage{esint}
\usepackage{algorithm}
\usepackage{algorithmic}
\newtheorem{prop}{Proposition}
\newtheorem{theo}{Theorem}
\newenvironment{dem}[1][Proof.]{\begin{trivlist}
\item[\hskip \labelsep {\bfseries #1}]}{\end{trivlist}}
\newcommand{\vu}{{\bf{u}}}
\newcommand{\vf}{{\bf{f}}}
\def\bcdot{{*}}

\begin{document}

\begin{frontmatter}

\title{\textbf{\large Compact Approximate Taylor methods for systems of conservation laws.}}

\author{H. Carrillo\fnref{myfootnote}} \author{C. Par\'es \fnref{myfootnote}}
\address{University of M\'alaga}
\fntext[myfootnote]{pares@uma.es, hugo.carrillo@uma.es}

\begin{abstract}

A new family of high order methods for systems of conservation laws are 
introduced: the Compact Approximate Taylor (CAT) methods.
These methods are based on centered $(2p + 1)$-point stencils where $p$ 
is an arbitrary integer.
We prove that the order of accuracy is $2p$ and that CAT methods are an 
extension of high-order Lax-Wendroff methods for linear problems.
Due to this, they are linearly $L^2$-stable under a $CFL-1$ condition. In 
order to prevent the spurious oscillations that appear close to 
discontinuities two shock-capturing techniques have been considered: a 
flux-limiter technique (FL-CAT methods) and WENO reconstruction for the 
first time derivative (WENO-CAT methods). We follow  \cite{ZORIO} in 
the second approach. A number of test cases are considered to compare 
these methods with other WENO-based schemes: the 
linear transport equation, Burgers equation,
and the 1D compressible Euler system are considered. Although CAT 
methods present an extra computational cost due to the local character, 
this extra cost is compensated by the fact that they still give good 
solutions with CFL values close to 1.

\end{abstract}

\end{frontmatter}

\section{Introduction}

Lax-Wendroff methods  for linear systems of conservation laws are based on Taylor expansions in time in which the time derivatives are transformed into spatial derivatives using the equations  \cite{Lax1}, \cite{LeVeque:2007:FDM:1355322},\cite{toro2009riemann}, \cite{GIDEON10.2307/2004918}. The spatial derivatives are then discretized by means of centered high-order differentiation formulas. This procedure allows to derive numerical methods of order $2p$, where $p$ is an arbitrary integer, using a centered  $(2 p + 1)$-point stencil that are $L^2$-stable (a review on these methods will be presented in Section 2).\\
This paper focuses on the extension of Lax-Wendroff methods to nonlinear systems of conservation laws. This problem is closely related to the design of numerical schemes based on the methods of lines in which the time discretization is performed by means of Taylor or Approximate Taylor methods. Many authors have focused on the design of this type of methods that can be an alternative to time discretizations based on SSP Runge-Kutta methods \cite{JIANXIAN1064827502412504} that lead to some stability restrictions for orders bigger than three.
The main difficulty to extend Lax-Wendroff methods to nonlinear problems come from the transformation of time derivatives into spatial derivatives using the equations. A first strategy to do this is given by the Cauchy-Kovalevskaya (CK) procedure. In \cite{JIANXIAN1064827502412504} this procedure has been used together with WENO reconstructions for  the spatial discretization. The main benefit compared to RK time discretizations is that only one WENO reconstruction is needed at each spatial cell per time step. On the other hand, the main drawback comes from the fact that the  CK procedure leads to expressions whose number of terms grow exponentially what implies high computational costs and difficult implementations. In \cite{ZORIO} an alternative has been presented based on an Approximate Taylor (AT) method in which the time derivatives are approximated using high-order centered differentiation formulas combined with Taylor approximations in time that are computed in a recursive way. The resulting method is easy to implement and shows a good performance. Nevertheless AT schemes are not proper generalizations of Lax-Wendroff methods: they have $(4p + 1)$-point stencils and worse linear stability properties than the original Lax-Wendroff methods. Nevertheless, they can be stabilized by using one WENO reconstruction  per spatial cell and time step, like in \cite{JIANXIAN1064827502412504} the resulting methods give good results under a $CFL-0.5$ condition typically.\\
In order to design numerical methods that are proper generalization of Lax-Wendroff methods,  a compact variant of the AT procedure is introduced here: first, the conservative expression of the high-order 
Lax-Wendoff methods is considered; then, the derivatives appearing in the expression of the numerical flux are computed using  $2p$-point differentiation formulas. This strategy lead to Compact Approximate 
Taylor (CAT) methods that have $2 p +1$-point stencil and order of accuracy $2p$. They reduce to the Lax-Wendroff method when applied to a linear systems and thus they are linearly $L^2$-stable under a $CFL-1$ condition. Nevertheless, the number of operations to perform a step is bigger than the original AT methods: this is due to the fact that the approximation of the time derivatives is local in the sense that they depend both on the point and on the stencil. However, this extra cost is compensated by better stability properties. As it happens with its linear counterpart, CAT methods lead to spurious oscillations near discontinuities. In order to  cure them, two shock-capturing techniques are considered here: a flux-limiter technique and the use of a WENO reconstruction per cell and time step, as in \cite{ZORIO} and  \cite{JIANXIAN1064827502412504}.\\
This paper is organized as follows: in Section 2 a review of high order Lax-Wendroff methods for the linear transport equation is presented, including the study of the order and the $L^2$-stability as well as the computation and properties of the coefficients. Section 3 is devoted to their extension to nonlinear problems: first the AT technique is recalled and then CAT methods are presented. We show that they reduce to Lax-Wendroff methods when applied to a linear problem and we analyze the order of accuracy. In Section 4 the techniques considered to cure the spurious oscillations near the discontinuities are presented. In Section 5 CAT methods are compared in a number of test cases with WENO-RK methods and AT methods. The linear transport equation, Burgers equation, and the 1D compressible Euler system are considered. Future developments and conclusions are drawn in Section 6.

\section{The high-order Lax-Wendroff method for linear problems}
Let us first consider the linear scalar equation:
\begin{equation}\label{edplinsc}
u_t+au_x=0.
\end{equation}
We consider the numerical method:
\begin{equation}\label{esquema_original}
u^{n+1}_i = u^n_i+\sum_{k=1}^{m} \frac{(-1)^k c^k}{k!} \sum_{j=-p}^{p} \delta^k_{p,j} \: u^n_{i+j},
\end{equation}
where $\{ x_i \}$ are the nodes of a uniform mesh of step $\Delta x$; $u_i^n$ is an approximation of the point value of the solution at $x_i$ at the time $n \Delta t$, where $\Delta t$ is the time step; $p \geq 1$ is a natural number;
$c= {a \Delta t }/{\Delta x}$; and $\delta^k_{p,j}$ are the coefficients of the centered interpolatory formula of numerical differentiation  based on a 
$(2p + 1)$-point  stencil:
\begin{equation} \label{F}
f^{(k)}(x_i)   \simeq   D^k_{p,i}(f, \Delta x) =  \frac{1}{\Delta x^k} \sum_{j=-p}^{p} \delta^k_{p,j} f(x_{i+j}),
\end{equation}
where $f^{(k)}$ represents the $k$-th derivative of a one-variable function $f$ and $f^{(0)} = f$. 
The expression of the numerical method is obtained by applying a Taylor expansion in time, and replacing time derivatives   by space derivatives through the identities
\begin{equation}\label{permder}
\partial^k_t u = (-1)^k a^k \partial^k_x u,\quad k=1, 2\dots
\end{equation}

\subsection{Formulas of numerical differentiation}
Besides \eqref{F} the following family of interpolatory formulas based on a $2p$-point stencil will be used in this work:
\begin{equation}\label{upwF}
f^{(k)}(x_i + q \Delta x) \simeq A^{k,q}_{p,i}(f, \Delta x) = \frac{1}{\Delta x^k} \sum_{j = -p + 1}^p \gamma^{k,q}_{p,j} f(x_{i+j}),
\end{equation}
i.e. $A^{k,q}_{p,i}(f, \Delta x)$ is the numerical differentiation formula that approximates the $k$-th derivative at the point $x_i + q \Delta x$ using the values of the 
function at the $2p$  points $x_{i-p+1}, \dots, x_{i+p}$. Observe that the coefficients, like in \eqref{F}, do not depend on $i$.

Given a variable $w$, the following  notation will be used:
\begin{eqnarray*}
 D^k_{p,i}(w_\bcdot, \Delta x)  & = &  \frac{1}{\Delta x^k} \sum_{j=-p}^{p} \delta^k_{p,j} w_{i+j},\\
 A^{k,q}_{p,i}(w_\bcdot, \Delta x)  & = & \frac{1}{\Delta x^k} \sum_{j = -p + 1}^p \gamma^{k,q}_{p,j} w_{i+j},
\end{eqnarray*}
to indicate that the formulas are applied to the approximations of $w$, $w_i$, and not to its exact point values $w(x_i)$. 
In cases where there are two or more indexes, the symbol $\bcdot$ will be used to indicate with respect to which the differentiation is applied. For instance:
\begin{eqnarray*}
\partial^k_x u(x_i + q \Delta x, t_n) & \simeq & A^{k,q}_{p,i} (u_\bcdot^n, \Delta x) =  \frac{1}{\Delta x^k} \sum_{j = -p + 1}^p \gamma^{k,q}_{p,j} u^n_{i+j}, \\
\partial^k_t u(x_i, t_n + q\Delta t) &\simeq& A^{k,q}_{p,n} (u_i^\bcdot, \Delta t) =  \frac{1}{\Delta t^k} \sum_{r = -p + 1}^p \gamma^{k,q}_{p,r} u^{n+r}_{i}.
\end{eqnarray*}
Using this notation, the algorithm \eqref{esquema_original} writes as follows:
\begin{equation}\label{lwndf}
u^{n+1}_i = u^n_i+\sum_{k=1}^{m} \frac{(-1)^k a^k  \Delta t^k}{k!} D^k_{p,i}(u^n_\bcdot, \Delta x).
\end{equation}

Let us discuss some properties of the coefficients of the numerical differentiation formulas \eqref{F} and \eqref{upwF} and some relations between them that will be used in that follows. Since the coefficients are independent
of $\Delta x$ and $i$, we can consider, without loss of generality, the case $i = 0$, $x_0 = 0$, $\Delta x = 1$:
\begin{equation} \label{ndf}
f^{(k)} (0)   \simeq    D^k_{p,0}(f,1) =  \sum_{j=-p}^{p} \delta^k_{p,j} f(j),
\end{equation}
\begin{equation}\label{upwndf}
f^{(k)}(q) \simeq  A^{k,q}_{p,0}(f,1) = \sum_{j = -p + 1}^p \gamma^{k,q}_{p,j} f(j).
\end{equation}

Since \eqref{ndf} is exact for polynomials of degree $\leq 2p$, by applying the formula to $x^s$, $s = 0, \dots, 2p$ at $x = 0$, we get that the coefficients have to satisfy the equalities
\begin{equation}\label{propbasc}
\sum_{j=-p}^{p} j^k \delta^k_{p,j}= k!, \quad\quad
\sum_{j=-p}^{p} j^s \delta^k_{p,j}= 0, \quad
s \neq k \quad,  \quad 0 \leq s,k \leq 2p.
\end{equation}
Analogously:
\begin{equation}\label{propbasc2}
\sum_{j=-p+1}^{p} j^k \gamma^{k,0}_{p,j}= k!, \quad\quad
\sum_{j=-p+1}^{p} j^s \gamma^{k,0}_{p,j}= 0, \quad
s \neq k \quad,  \quad 0 \leq s,k \leq 2p-1.
\end{equation}
\begin{equation}\label{propbasc3}
\sum_{j=-p+1}^{p}  \gamma^{k,q}_{p,j}=\left\{ \begin{array}{ll}  1 & \text{if $k = 1$,} \\ 0 & \text{otherwise.}  \end{array} \right.
\end{equation}

As it is well known, the coefficients $\delta^k_{p, j}$ are related to the Lagrange basis polynomials
\begin{equation}\label{lagragianfunction}
  F_{p,j}(x)= \prod_{r=-p,r\neq j}^{p} \frac{ (x-r)}{(j-r)}, \quad -p\leq j \leq p,
\end{equation}
through the equalities:
\begin{equation}\label{deltas}
\delta^k_{p,j}  = F^{(k)}_{p,j}(0), 
\end{equation}
which allow us to write the Taylor expansion of $F_{p,j}$ centered at $x = 0$  as follows: 
\begin{equation}\label{taylorform}
F_{p,j}(x) = \sum_{k=0}^{2p} \frac{\delta^k_{p,j} }{k!}x^k.
\end{equation}

\begin{prop}
The coefficients $\delta^k_{p,j}$ of the formula (\ref{ndf}), satisfy:
\begin{eqnarray}
& &  \delta^k_{p,j}= (-1)^k \delta^k_{p,-j}; \label{p1}\\
& &  \delta^k_{p,0}= 0 \text{ if  $k$ is odd}; \label{p2}\\
& & \sum_{j=-p}^{p}\delta^k_{p,j} \, j^{(2p+1)}= 0\text{ if $k$ is even}; \label{p3}\\
& &   \sum_{j=-p}^{p}\delta^k_{p,j} \, j^{(2p+2)}= 0 \text{ if  $k$ is odd}. \label{p4}
\end{eqnarray}
\end{prop}

\begin{dem}
\eqref{p1} is deduced from the equality  $$ F_{p,-j}(x) = F_{p,j}(-x);$$. Using \eqref{p1} we get  \eqref{p2}. \eqref{p3} and \eqref{p4} are deduced from \eqref{p1}  and \eqref{p2}.
\end{dem}
\hfill$\Box$

\begin{prop}\label{prop:deltad}
For $k \geq 1$ the following relations hold:
\begin{eqnarray}
& &  \delta^k_{p,p} = \gamma^{k-1, 1/2}_{p,p}; \\
& &  \delta^k_{p,j} =  \gamma^{k-1, 1/2}_{p,j} - \gamma^{k-1, 1/2}_{p, j+1}, \quad j = -p+1, \dots, p-1;\\
 & & \delta^k_{p,-p} = - \gamma^{k-1, 1/2}_{p,-p+1}. 
\end{eqnarray}
\end{prop}

\begin{dem}
Let us consider the formulas
\begin{equation}\label{f1/2}
f^{(k-1)}(1/2) \simeq  A^{k-1, 1/2}_{p,0}(f,1) = \sum_{j = -p + 1}^p \gamma^{k-1}_{p,j} f(j),
\end{equation}
\begin{equation}\label{f-1/2}
f^{(k-1)}(-1/2) \simeq  A^{k-1, 1/2}_{p,-1}(f,1) = \sum_{j = -p + 1}^p \gamma^{k-1}_{p,j} f(j-1),
\end{equation}
that are exact for polynomials of degree $\leq 2p-1$.
Let us consider now the formula
\begin{equation}\label{ndf1}
f^{(k)}(0) \simeq A^{k-1, 1/2}_{p,0}(f,1) -  A^{k-1, 1/2}_{p,-1}(f,1).
\end{equation}
If $f$ is a polynomial of degree $2p $, then \eqref{f1/2} and \eqref{f-1/2} are exact for $f$, furthermore
$$
 A^{k-1, 1/2}_{p,0}(f,1) -  A^{k-1, 1/2}_{p,-1}(f,1) = f^{(k-1)} (1/2) -  f^{(k-1)}(-1/2) = f^{(k)}(0),
$$
where we have used that the formula
$$
g'(0) \simeq g(1/2) - g(-1/2),
$$
is exact for polynomials of degree 1.  Therefore, \eqref{ndf1} coincide with \eqref{ndf}. The proof is finished by writing \eqref{ndf1} in the form
\begin{align*}
f^{(k)}(0)  \simeq & \gamma^{k-1, 1/2}_{p,p} f(p) + (\gamma^{k-1, 1/2}_{p, p-1} - \gamma^{k-1, 1/2}_{p,p}) f(p-1) +  \dots \\
  & + (\gamma^{k-1, 1/2}_{p, -p+1} - \gamma^{k-1, 1/2}_{p, -p+2}) f(-p+1) - \gamma^{k-1, 1/2}_{p, -p+1}f(-p),
\end{align*}
and matching the weights. 
\end{dem}
\hfill$\Box$

\begin{prop} \label{prop3}
Given $1 \leq k \leq 2p-1$, $0 \leq s \leq k$:
\begin{equation}\label{CsCk}
\sum_{j=-p +1}^p \gamma^{s,q}_{p,j} \gamma^{k-s,j}_{p,l} = \gamma^{k,q}_{p,l}, \quad l= -p+1, \dots, p.
\end{equation}
\end{prop}

\begin{dem}
The proof is similar to the one of the preceding in  Proposition \ref{prop:deltad}: consider the formula
$$
f^{(k)}(q) \simeq \sum_{j=-p+1}^p \gamma^{s,q}_{p,j} f^{(k-s)}_j,
$$
with
$$
 f^{(k-s)}_j = \sum_{l = -p +1}^p \gamma^{k-s,j}_{p,l}f(l);
$$
check that it is exact for polynomials of degree $2p-1$; write it in the form:
$$
f^{(k)}(q) \simeq \sum_{l = -p+1}^p \left(\sum_{j=-p +1}^p \gamma^{s,q}_{p,j} \gamma^{k-s,j}_{p,l}  \right) f(l);
$$
and match its weights with those of \eqref{upwndf}.
\end{dem}
\hfill$\Box$

\subsection{Conservative form of (2)}
From the proof of Proposition \ref{prop:deltad} we deduce an alternative form for \eqref{F}:
\begin{equation}\label{F2}
f^{(k)}(x_i) \simeq  \frac{1}{\Delta x} \left(A^{k-1, 1/2}_{p,i}(f, \Delta x) -  A^{k-1, 1/2}_{p,i-1}(f, \Delta x)\right).
\end{equation}
Using this form in \eqref{lwndf},  the numerical 
method \eqref{esquema_original} can be written as:
\begin{equation}\label{cons}
u_i^{n+1} = u_i^n + \frac{\Delta t}{\Delta x}\left( F^p_{i-1/2} - F^p_{i+1/2}\right),
\end{equation}
with
\begin{equation}\label{flujoconservativo}
F^p_{i+1/2} =  \sum_{k=1}^{2p} (-1)^{k-1} \frac{a^k \Delta t^{k-1}}{k!} A^{k-1, 1/2}_{p,i}(u^n_\bcdot, \Delta x) .
\end{equation}
Using  \eqref{propbasc3} it is straightforward to verify that $F^p_{i+1/2}$ is a consistent numerical flux, what proves that \eqref{esquema_original} is a conservative method.

\subsection{Computation of the coefficients: an iterative algorithm}
Notice that \eqref{propbasc} constitutes a $(2p+1)\times(2p+1)$ linear system with a Vandermonde matrix that can be used to compute   $\delta^k_{p,i}$ . Nevertheless, as it is well-known, this system is  ill-conditioned, so that it is recommendable to compute them by using an alternative algorithm: we adapt the recursive algorithm proposed in \cite{FORNBERG2008770}. 
The following notation is adopted:
\begin{equation*}
\delta^k_{p,j} = 0 \mbox{ if } k > 2p \mbox{ or }   k< 0.
\end{equation*}
Let us derive some recurrence formulas to compute the coefficients:
\begin{enumerate}
\item   $\delta^k_{p,j} $ for  $j=0,..,p-1$.\\
From (\ref{lagragianfunction}), we obtain
\begin{align}\label{ff2}
F_{p,j}(x)&= \frac{x^2-p^2}{j^2-p^2} F_{p-1,j}(x).
 \end{align}
Using then the Taylor expansions (\ref{taylorform}) in (\ref{ff2}) we get
\begin{equation}\label{deltacentrados}
\delta^k_{p,j}= \frac{1}{p^2-k^2} \left[ p^2 \delta^{k}_{p-1,j} - k(k-1)\delta^{k-2}_{p-1,j} \right],
\end{equation}

\item  $\delta^k_{p,j}$ with  $j=p$. \\
Substituting $\textit{j=p}$ in (\ref{lagragianfunction}), we get
\begin{equation}\label{lateral}
F_{p,p}(x) = \frac{1}{(2p)(2p-1)}(x^2+x-p(p-1))F_{p-1,p-1}(x),
 \end{equation}
and, using (\ref{taylorform}), we obtain:
\begin{equation}\label{deltafinal}
\delta^k_{p,p} = \frac{1}{2p(2p-1)} \left[ k(k-1)\delta^{k-2}_{p-1,p-1} + k \delta^{k-1}_{p-1,p-1}-p(p-1)\delta^{k}_{p-1,p-1} \right].
\end{equation}
\item $\delta^k_{p,j}$ for ${j=-p, \dots, -1}$.  $\eqref{p1}$ is used. 
\end{enumerate}
The algorithm is computed only once in increasing order of $p$. The coefficients $\gamma^{k,q}_{p,j}$ are computed using the algorithms described in \cite{FORNBERG2008770},\cite{Fornberg2} and $\gamma^{k,1/2}_{p,j}$ is obtained from $\delta^{k+1}_{p,j}$.

\subsection{Order of accuracy}

\begin{prop} \label{Proporderndf} The formula of numerical differentiation (\ref{F}) has order of accuracy $\alpha_k-k$, with, 
\begin{equation*}\label{proposition2}
\alpha_k = \left\{
\begin{array}{ll}
2p+1 & \mbox{if $ k$ is odd,}\\
2p+2 & \mbox{if $k$ is even.}
\end{array}\right.\\
\end{equation*}
\end{prop}
\begin{dem}
Let $f$ be a function of class $C^{\alpha_k +1}$. Applying Taylor expansions  and properties \eqref{propbasc} and \eqref{p3}, we obtain:
\begin{equation}
\frac{1}{\Delta x^{k}}\sum_{j=-p}^{p} \delta^k_{p,j} f(x_{i+j}) =f^{(k)}(x_i)+ \varphi_k \frac{\Delta x^{\alpha_k-k}}{\alpha_k!} f^{(\alpha_k)} (x_i)+\mathcal{O}(\Delta x^{\alpha_k-k+1}),
\end{equation}
where
\begin{equation}\label{gamma}
\varphi_k = \sum_{j=-p}^{p} \delta^k_{p,j} \, j^{\alpha_k}.
\end{equation}
\hfill$\Box$
\end{dem}
Table \ref{tabla1} shows the order of \eqref{F} for different values of $p$ and $k$.
\renewcommand{\tablename}{Table}
\begin{table}
\centering
\begin{tabular}{|l||*{4}{c|}}\hline
\backslashbox{k}{p}
&\makebox[2em]{1}&\makebox[3em]{2}&\makebox[3em]{3}
&\makebox[3em]{4}\\\hline\hline
1 &2&4&6&8\\\hline
2 &2&4&6&8\\\hline
3 &&2&4&6\\\hline
4 &&2&4&6\\\hline
5 &&&2&4\\\hline
6 &&&2&4\\\hline
7 &&&&2\\\hline
8 &&&&2\\\hline
\end{tabular}
\caption{Order of the formula (\ref{F}).}
\label{tabla1}
\end{table}

\begin{prop} The discretization error of the numerical method \eqref{esquema_original} is of order $O(\Delta t^{m+1} + \Delta x^{2p +1})$.
\end{prop}

\begin{dem} Let $u$ be a smooth enough solution of \eqref{edplinsc}. Using Proposition \ref{Proporderndf} we obtain
\begin{eqnarray*}
& & u(x_i, t_{n+1}) - u(x_i, t_n) - \sum_{k=1}^{m} \frac{(-1)^k c^k}{k!} \sum_{j=-p}^{p} \delta^k_{p,j} \: u(x_{i+j}, t_n) \\
& & \quad =  u(x_i, t_{n+1}) - u(x_i, t_n) \\
& & \qquad - \sum_{k=1}^{m} \frac{(-1)^k a^k \Delta t^k}{k!} \Bigl( \partial_x^ku(x_i, t_n) + 
\varphi_k \frac{\Delta x^{\alpha_k-k}}{\alpha_k!} \partial_x^{\alpha_k} u(x_i,t_n) + \mathcal{O}(\Delta x^{\alpha_k-k+1})\Bigr)\\
& & \quad =  u(x_i, t_{n+1}) - u(x_i, t_n)  - \sum_{k=1}^{m} \frac{\Delta t^k}{k!} \partial_t^k u(x_i, t_n) \\
& & \qquad - \sum_{k=1}^m \varphi_k \frac{(-1)^k c^k }{k!\alpha_k!} \Delta x^{\alpha_k} \, \partial^{\alpha_k}_xu(x_i, t_n) + O(\Delta x^{\alpha_k + 1})\\
& & \quad = \frac{1}{(m+1)!} \partial_t^{m+1} u(x_i, t_n) \Delta t^{m+1} \\
& & \qquad +\left(  \sum_{k=0}^{p-1} \frac{ \varphi_{2k+1} c^{2k+1} } {(2p+1)!(2k+1)!}\right)    \partial^{2p +1}_x  u(x_i, t_n)  \Delta x^{2p+1}+ O(\Delta t^{m+2} 
+ \Delta x^{2p + 2}),
\end{eqnarray*}
where \eqref{permder} has been used. 
\end{dem}
\hfill$\Box$

As a consequence, the order of accuracy of \eqref{esquema_original} is  $\min(m, 2p)$. Therefore, the optimal combination of these parameters is  $m = 2p$. 
From now on, we shall assume that this relation holds.

\subsection{Modified equation and stability}
Taking into account that $m = 2p$ and \eqref{p4}, the local discretization error is as follows:
\begin{eqnarray*}
& & u(x_i, t_{n+1}) - u(x_i, t_n) - \sum_{k=1}^{m} \frac{(-1)^k c^k}{k!} \sum_{j=-p}^{p} \delta^k_{p,j} \: u(x_{i+j}, t_n) \\
& & \quad = \frac{1}{(2p+1)!} \partial_t^{2p+1} u(x_i, t_n) \Delta t^{2p+1} + \frac{1}{(2p+2)!} \partial_t^{2p+2} u(x_i, t_n) \Delta t^{2p+2}\\
& & \qquad +\left(  \sum_{k=0}^{p-1} \frac{ \varphi_{2k+1} c^{2k+1} } {(2p+1)!(2k+1)!}  \right)  \partial^{2p +1}_x  u(x_i, t_n)  \Delta x^{2p+1} \\
& & \qquad - \left(  \sum_{k=1}^{p} \frac{ \varphi_{2k} c^{2k} } {(2p+2)!(2k)!}   \right)  \partial^{2p +2}_x  u(x_i, t_n) \Delta x^{2p+2} + O(\Delta x^{2p + 3}).
\end{eqnarray*}
Using \eqref{gamma} and \eqref{taylorform} we get:
\begin{eqnarray*}
\sum_{k=0}^{p-1} \frac{\varphi_{2k+1} c^{2k+1}}{(2k+1)!}  & =  & \sum_{j=-p}^{p} \left( \sum_{k=0}^{p-1} \frac{\delta^{2k+1}_{p,j} c^{2k+1}}{(2k+1)!} \right) j^{2p+1}\\
&  = &  \frac{1}{2}  \sum_{j=-p}^{p} \left( \sum_{l=1}^{2p} \left( \frac{\delta^{l}_{p,j}}{l!}-\frac{\delta^{l}_{p,-j}}{l!} \right)  c^{l} \right) j^{2p+1}\\
 & =  & \frac{1}{2}  \sum_{j=-p}^{p} \left( F_{p,j}(c)-F_{p,-j}(c)\right) j^{2p+1}\\
&  = &  \frac{1}{2} \left[\sum_{j=-p}^{p}  F_{p,j}(c) j^{2p+1} - \sum_{j=-p}^{p} F_{p,j}(-c) j^{2p+1} \right]\\
& =  &  \frac{1}{2} \left[q(c) - q(-c)\right],\\
\end{eqnarray*}
where $q(c)$ is the polynomial of degree $\leq 2p$ that interpolates the points  
$$\{(-p,(-p)^{2p+1}),\dots,(0,0),\dots,(p,p^{2p+1}) \} .$$
Since $q$ is clearly an odd function, we finally obtain:
\begin{equation}\label{q}
\sum_{k=0}^{p-1} \frac{\varphi_{2k+1} c^{2k+1}}{(2k+1)!}  = q(c).
\end{equation}
Reasoning in a similar way, we get:
\begin{equation}\label{r}
\sum_{k=1}^{p} \frac{ \varphi_{2k} c^{2k} } {(2k)!}   = r(c),
\end{equation}
where $r$   is the polynomial of degree $\leq 2p$ that interpolates the points 
$$\{(-p,(-p)^{2p+2}),...,(0,0),...,(p,p^{2p+2}) \} .$$
Using now \eqref{permder}, \eqref{q}, and \eqref{r}, the local discretization error can be written as follows:
\begin{eqnarray*}
& & u(x_i, t_{n+1}) - u(x_i, t_n) - \sum_{k=1}^{m} \frac{(-1)^k c^k}{k!} \sum_{j=-p}^{p} \delta^k_{p,j} \: u(x_{i+j}, t_n) \\
& & \qquad =  \frac{   h_1(c) }{(2p+1)!}  \partial_x^{2p+1}u(x_i,t_n) \Delta x^{2p+1}  -\frac{ h_2(c)}{(2p+2)!}   \partial^{2p+2}_x u(x_i,t_n) \Delta x^{2p+2} +\mathcal{O} (\Delta^{2p+3}) ,
\end{eqnarray*}
with
\begin{align}
&h_1(c)=q(c)-c^{2p+1},\\
&h_2(c)=r(c)-c^{2p+2}. 
\end{align}

Therefore, the numerical method solves with order  $\mathcal{O} (\Delta^{2p+2})$ the following  modified equation
\begin{equation}\label{eqmodLW}
u_t + a u_x =  \mu_1  \partial_x^{2p+1}u    - \mu_2     \partial^{2p+2}_x u,
\end{equation}
where
\begin{equation}\label{eqmodLW2}
\mu_1 =   \frac{   h_1(c) }{(2p+1)!\Delta t}\Delta x^{2p+1}, \qquad \mu_2 = \frac{ h_2(c)}{(2p+2)!\Delta t}  \Delta x^{2p+2} .
\end{equation}

In order to study the stability, let  us look for an elementary solution $u(x,t)$ of \eqref{eqmodLW2} of the form 
$$
u(x,t)  = e^{\alpha t}\cdot e^{i kx},
$$ 
where $\alpha$ is complex number. The following equality has to be satisfied:
$$
\alpha u  +  ik a  u =   \mu_1 (-1)^{p} i k^{2p+1} u  + \mu_2 (-1)^{p}k^{2p+2} u.
$$
Therefore:
$$\alpha =  - \mu_2 (-1)^{p+1}k^{2p+2} - (k a - \mu_1 (-1)^{p+1}k^{2p+1} ) i.$$
The numerical method is thus expected to be stable if the real part is negative, i.e.
$$
 \mu_2 (-1)^{p} \leq 0, 
 $$
or, equivalently
\begin{equation}\label{stabcond}
h_2(c) (-1)^{p} \leq 0. 
\end{equation}

\begin{figure}[!ht]
	\small
	\setlength{\unitlength}{1mm}
	\centering	
	\begin{picture}(70,65)
	\put(-17,3){\makebox(100,50)[c]{
	\includegraphics[width=19cm,height=7cm]{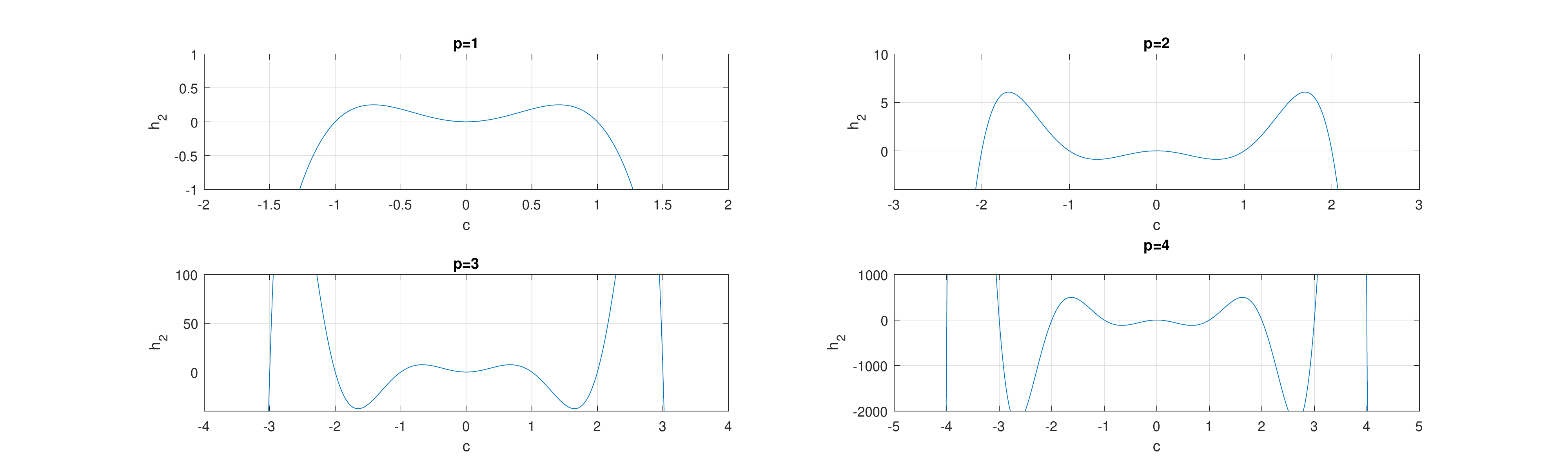}}}
	\end{picture}
	\vspace{0 cm}
	\caption{Function $h_2(c)$ for $p=1,..,4$.}
	\label{figura_ex_2}
	\end{figure}

$h_2$ is a pair polynomial of degree $2p + 2$  such that
$$
\lim_{c \rightarrow \pm \infty} h_2(c) = -\infty.
$$
Moreover, $0$ is a double root of $h_2$ and $\pm 1, \dots, \pm p$ are single roots. Analyzing the change of signs of $h_2$,  it can be easily checked that
\begin{eqnarray*}
&& h_2(c)\leq 0 ,\quad  \forall c\in [0,1]  \quad \text{if $p$ even}, \\
&& h_2(c)\geq 0, \quad \forall c\in [0,1]  \quad \text{if $p$ odd},
\end{eqnarray*}
and thus \eqref{stabcond} is satisfied if $c \in [0,1]$ (see Figure \ref{figura_ex_2}). 
This argument shows that the stability of the method  under the standard CFL condition $c \leq 1$ can be formalized: extended details and a formal proof can be found in \cite{JIEQUANLI2013610}. A study of the modified equation for the Lax Wendroff method can also be found  in \cite{WARMING1974159}. 

\section{Extension to nonlinear problems}

\subsection{Approximate Taylor method}

Following  \cite{ZORIO}, instead of using the Cauchy-Kovaleskaya process to extend \eqref{cons}-\eqref{flujoconservativo} to nonlinear problems
\begin{equation}\label{scl}
u_t+f(u)_x=0,
\end{equation}
we use the equalities
\begin{equation}
\partial_t^k u = -\partial_x\partial_t^{k-1}f(u).
\end{equation}
To derive the expression of the numerical method, let us suppose that approximations
$$
\tilde{f}^{(k-1)}_{i}  \cong \partial_t^{k-1}f(u)(x_{i}, t_n),
$$
are available. Then,
$$
\partial_t^k u(x_i, t_n) \cong  \tilde u^{(k)}_i =  - D^1_{p_{k-1},i} (\tilde{f}^{(k-1)}_\bcdot, \Delta x) = - \frac{1}{\Delta x}\ \sum_{j=-p_{k-1}}^{p_{k-1}} \delta^1_{p_{k-1},j}  \tilde{f}^{(k-1)}_{i+j} ,
$$
being

\begin{equation}\label{p_k}
p_{k} = \lceil (p-k/2) \rceil,
\end{equation}
where, $\lceil \cdot \rceil $ denotes the ceiling function.\\
Using these approximations to approximate the Taylor expansion, we obtain the method
\begin{equation}\label{atm}
u_i^{n+1} = u_i^n +\sum_{k= 0}^{2p} \frac{\Delta t ^k}{k! }\  \tilde u^{(k)}_i.
\end{equation}
Equivalently, using \eqref{F2}, the numerical method can be written in conservative form \eqref{cons} with numerical flux
\begin{equation}
\label{nlcase1}
F^p_{i+1/2}  = \sum_{k=1}^{2p} \frac{\Delta t^{k-1}}{k!}A^{0, 1/2}_{p_{k-1},i}(\tilde{f}^{(k-1)}_\bcdot, \Delta x),   
\end{equation}
being
\begin{equation}\label{nlcase2}
A^{0,1/2}_{p_{k-1},i}(\tilde{f}^{(k-1)}_\bcdot, \Delta x) = \sum_{j=-p_{k-1}+1}^{p_{k-1}} \gamma^{0,1/2}_{p_{k-1},j}  \tilde{f}^{(k-1)}_{i+j} .
\end{equation}

Now, to compute the approximations $\tilde{f}^{(k-1)}_{i}$, new Taylor expansions in time are used recursively as follows:

\begin{itemize}

\item  Compute $\tilde{f}^{(0)}_{i}=f(u^n_{i})$.

\begin{itemize}

\item  For $k = 1 \dots 2p$:

\begin{itemize}

\item Compute
$$ \tilde{u}  ^{(k)}_{i} = - D^1_{p_{k-1},i}(\tilde{f}^{(k-1)}_\bcdot, \Delta x).$$

\item Compute
$$
\tilde{f}^{k,n+r}_{i} = f \left(  u^n_{i} + \sum_{l=1}^{k} \frac{(r \Delta t)^l}{l!} \tilde{u}^{(l)}_{i} \right), \quad r = -p, \dots, p.
$$

\item Compute
$$
\tilde{f}^{(k)}_{i} =  D^k_{p,n}(\tilde{f}^{k, \bcdot}_i, \Delta t),
$$
where
$$
D^k_{p,n}(\tilde{f}^{k,\bcdot}_i, \Delta t) = \frac{1}{\Delta t^k} \sum_{r= -p}^p \delta^k_{p,r} \tilde{f}^{k,n+r}_{i}.
$$

\end{itemize}

\end{itemize}

\end{itemize}

Observe that Taylor expansions are used to approximate $f(u(x_i, t_n + r \Delta t))$ and  once these approximations have been computed, the centered formula of numerical differentiation \eqref{F} is used to approximate the temporal derivatives. 

This method is order $2p$, but it is not a generalization of \eqref{esquema_original} in the sense that this latter method is not recovered if $f(u) = au$. To see this, consider $p = 1$ and $f(u) = au$: it can be easily checked that \eqref{atm} writes as follows 
\begin{equation}\label{5plwm}
u_i^{n+1} = u_i^n - \frac{c}{2}(u_{i+1}^n - u_{i-1}^n)  - \frac{c^2}{8}(u_{i+2}^n - 2 u_i^n  + u_{i-2}^n),
\end{equation}
which is different from the standard Lax-Wendroff method: \eqref{atm} is a $(4p + 1)$-point method whose stability properties are worse than those of the standard Lax-Wendroff method. (see \cite{LeVeque:2007:FDM:1355322}). 

\subsection{Compact Approximate Taylor method}\label{sec:CAT}
In order to prevent the increase of the stencil observed for Approximate Taylor methods, we consider a modification based on the conservative form of the method.
The numerical flux $F^p_{i+1/2}$ will be computed using only the approximations 
\begin{equation} \label{fluxstencil}
u^n_{i -p +1}, \dots, u^n_{i+ p}, 
\end{equation}
so that the values used to update $u_i^{n+1}$ are only those of the centered $(2 p + 1)$-point stencil, like in  the linear case. In fact, we will show that this modification is a proper generalization of the Lax-Wendroff method for linear problems. 

In order to be able to compute the numerical fluxes using only \eqref{fluxstencil}, for every $i$ we will compute \textit{local} approximations
of
$$
\partial_t^{k-1} f(u(x_{i-p+1}, t^n), \dots, \partial_t^{k-1} f(u(x_{i +p}, t^n), 
$$
that will be represented by
$$
\tilde{f}^{(k-1)}_{i, j}  \cong \partial_t^{k-1}f(u)(x_{i + j}, t_n), \quad j=-p+1, \dots, p.
$$
These approximations are local  in the sense that 
$i_1 + j_1 = i_2 + j_2$, does not imply that $\tilde{f}^{(k-1)}_{i_1, j_1}  = \tilde{f}^{(k-1)}_{i_2, j_2}  $. 
Once these approximations have been computed, the numerical flux is given by
\begin{equation}\label{cat1}
F^p_{i+1/2}  = \sum_{k=1}^{2p} \frac{\Delta t^{k-1}}{k!}A^{0, 1/2}_{p,0}(\tilde{f}_{i, \bcdot}^{(k-1)}, \Delta x),   
\end{equation}
with 
\begin{equation}\label{cat2}
A^{0, 1/2}_{p,0}(\tilde{f}_{i, \bcdot}^{(k-1)}, \Delta x) = \sum_{j=-p+1}^{p} \gamma^{0,1/2}_{p,j}  \tilde{f}^{(k-1)}_{i,j} .
\end{equation}

Now, given $i$,  to compute the approximations $\tilde{f}^{(k-1)}_{i, j}$, new Taylor expansions in time are used recursively as follows:

\begin{itemize}

\item  Compute $\tilde{f}^{(0)}_{i,j}=f(u^n_{i+j}), \quad  j = -p+1, \dots, p.$

\begin{itemize}

\item  For $k = 1 \dots 2p$:

\begin{itemize}

\item Compute
$$ \tilde{u}  ^{(k)}_{i,j} = - A^{1,j}_{p,0}(\tilde{f}^{(k-1)}_{i, \bcdot}, \Delta x),$$
where 
$$  A^{1,j}_{p,0}(\tilde{f}^{(k-1)}_{i,\bcdot}, \Delta x) = \frac{1}{\Delta x} \sum_{r=-p +1}^p \gamma^{1,j}_{p,r} \tilde{f}^{(k-1)}_{i, r}.$$ 

\item Compute
$$
\tilde{f}^{k,n+r}_{i,j} = f \left(  u^n_{i+j} + \sum_{l=1}^{k} \frac{(r \Delta t)^l}{l!} \tilde{u}^{(l)}_{i,j} \right), \quad  j, r = -p+1, \dots, p.
$$

\item Compute
$$
\tilde{f}^{(k)}_{i,j} =   A^{k,0}_{p,n}( \tilde{f}^{k, \bcdot}_{i,j}, \Delta t),\quad  j = -p+1, \dots, p.
$$
with
$$
 A^{k,0}_{p,n}( \tilde{f}^{k, \bcdot }_{i,j}, \Delta t) = \frac{1}{\Delta t^{k}} \sum_{r=-p+1}^{p} \gamma_{p,r}^{k,0}\tilde{f}^{k,n+r}_{i,j} .
$$
\end{itemize}

\end{itemize}

Notice that, unlike the Approximate Taylor methods (in which all the derivatives were approximated using the centered $(2p +1)$-point formula), in this algorithm the 
stencil $x_{i-p+1}, \dots, x_{i+p}$ is used for the space derivatives and the  stencil
$t_{n-p+1}, \dots, t_{n+p}$ for the time derivative. 

\end{itemize}

\begin{theo} The compact approximate Taylor method reduces to \eqref{esquema_original} when $f(u) = au$.
\end{theo}

\begin{dem}

For $k >1$ we have:
\begin{eqnarray*}
\tilde{f}^{(k-1)}_{i,j} & =  &  \frac{1}{\Delta t^{k-1}} \sum_{r=-p+1}^{p} \gamma_{p,r}^{k-1,0} \tilde{f}^{k-1,n+r}_{i,j}\\
& = &   \frac{a}{\Delta t^{k-1}} \sum_{r=-p+1}^{p} \gamma_{p,r}^{k-1,0} \left(   u^n_{i+j} +  \sum_{l=1}^{k-1} \frac{(r \Delta t)^l}{l!} \tilde{u}^{(l)}_{i,j} \right)\\
& = &    \frac{a}{\Delta t^{k-1}} \left(  \left( \sum_{r=-p+1}^{p} \gamma_{p,r}^{k-1,0} \right) u^n_{i+j} +  \sum_{l=1}^{k-1} \frac{\Delta t^l}{l!}
 \left(  \sum_{r=-p+1}^{p} \gamma_{p,r}^{k-1,0}  r^l \right)  \tilde{u}^{(l)}_{i,j} \right)\\
& = & a  \tilde{u}^{(k-1)}_{i,j}, 
\end{eqnarray*}
where \eqref{propbasc2} has been used. 
On the other hand:
\begin{eqnarray*}
\tilde{u}  ^{(k)}_{i,j} &  = & - \frac{1}{\Delta x}  \sum_{r=-p+1}^{p} \gamma_{p,r}^{1,j} \, \tilde{f}^{(k-1)}_{i,r} \\
 &  = & - \frac{a}{\Delta x}  \sum_{r=-p+1}^{p} \gamma_{p,r}^{1,j}  \tilde{u}^{(k-1)}_{i,r} \\
& = & \frac{a^2}{\Delta x^2} \sum_{r=-p+1}^{p} \gamma_{p,r}^{1,j} \sum_{s= -p+1}^p \gamma_{p,s}^{1,r} \tilde{u}^{(k-2)}_{i,s}\\
& = & \frac{a^2}{\Delta x^2}   \sum_{s= -p+1}^p \left(    \sum_{r=-p+1}^{p} \gamma_{p,r}^{1,j} \gamma_{p,s}^{1,r} \right) \tilde{u}^{(k-2)}_{i,s}\\
& = &  \frac{a^2}{\Delta x^2}   \sum_{s= -p+1}^p  \gamma_{p,s}^{2,j}  \tilde{u}^{(k-2)}_{i,s},
\end{eqnarray*}
where \eqref{CsCk} has been used. By recurrence:
\begin{equation}
\tilde{u}  ^{(k)}_{i,j} =\frac{(-1)^{k}a^{k}}{\Delta x^{k}}  \sum_{r=-p+1}^{p} \gamma_{p,r}^{k,j} u^n_{i+r} .
\end{equation}
Next,
\begin{eqnarray*}
A^{0,1/2}_{p,0}(\tilde{f}^{(k-1)}_{i, \bcdot}, \Delta x)  & = & \frac{1}{\Delta x}  \sum_{j=-p+1}^{p} \gamma_{p,j}^{0,1/2}  \tilde{f}^{(k-1)}_{i,j} \\
 & = & \frac{a}{\Delta x}  \sum_{j=-p+1}^{p} \gamma_{p,j}^{0,1/2}  \tilde{u}^{(k-1)}_{i,j} \\
& = &  (-1)^{k-1} \frac{a^{k}}{\Delta x^{k}}\sum_{j=-p+1}^{p} \gamma_{p,j}^{0,1/2}   \sum_{r=-p+1}^{p} \gamma_{p,r}^{k-1,j} u^n_{i+r}\\
& = & (-1)^{k-1} \frac{a^{k}}{\Delta x^{k}} \sum_{r=-p+1}^{p} \left(    \sum_{j=-p+1}^{p} \gamma_{p,j}^{0,1/2} \gamma_{p,r}^{k-1,j}   \right)  u^n_{i+r}\\
& = & (-1)^{k-1} \frac{a^{k}}{\Delta x^{k}} \sum_{r=-p+1}^{p}  \gamma_{p,j}^{k-1,1/2}  u^n_{i+r}\\
& = & (-1)^{k-1} a^{k} A^{k-1, 1/2}_{p,i}(u^n_\bcdot, \Delta x),
\end{eqnarray*}
where \eqref{CsCk} has been used. 
Finally,
\begin{eqnarray*}
F^p_{i+1/2}  & = & \sum_{k=1}^{2p} \frac{\Delta t^{k-1}} {k!} A^{0, 1/2}_{p,0}(\tilde{f}^{(k-1)}_{i, \bcdot}, \Delta x) \,\\
& = & \sum_{k=1}^{2p} (-1)^{k-1} \frac{a^k \Delta t^{k-1}}{k!} A^{k-1, 1/2}_{p,i}(u^n_\bcdot, \Delta x),  
\end{eqnarray*}
what is the numerical flux \eqref{flujoconservativo} corresponding to \eqref{esquema_original}, as we wanted to prove. 

\end{dem}
\hfill$\Box$

As a consequence, we obtain that the compact approximate Taylor method is linearly stable (in the $L^2$ sense) under the usual CFL condition
\begin{equation}\label{CFL}
\max_i(|f'(u_i)|) \frac{\Delta t}{\Delta x} \leq 1.
\end{equation}

\begin{theo} The compact approximate Taylor method is order $2p$.
\end{theo}
\begin{dem}

Let us perform  a step of the method starting from  the point values at time $t^n$, $u(x_i, t_n)$, of a smooth enough exact solution. We assume that $\Delta t/\Delta x$ remains constant.

First we have:
$$ \tilde{u}  ^{(1)}_{i,j} = - A^{1,j}_{p,0}(\tilde{f}^{(0)}_{i, \bcdot}, \Delta x) = -\partial_x f(u) (x_{i+j}, t_n) + O(\Delta x^{2p-1}) = \partial_t u(x_{i+j}, t_n) +  O(\Delta x^{2p-1}) . $$
Next:
$$
\tilde{f}_{i,j}^{1,n+r} = f( u(x_{i+j}, t_n) +  \tilde{u}  ^{(1)}_{i,j} r \Delta t ) = f(P^1_{i,j} (r \Delta t) ) + O(\Delta x^{2p}) ,
$$
where
$$
P^1_{i,j}(s) = u(x_{i+j}, t_n)  + s \partial_t u(x_{i+j}, t_n) ,
$$
is the first order Taylor polynomial in time of $u$ in $(x_{i+j}, t_n)$. 
Then
\begin{eqnarray*}
\tilde{f}^{(1)}_{i,j}  & = &   A^{1,0}_{p,n}( \tilde{f}^{k, \bcdot}_{i,j}, \Delta t)\\
& = & \frac{1}{\Delta t} \sum_{r= -p + 1}^p \gamma_{p,j}^{1,0} \tilde{f}^{1, n+r}_{i,j}\\
& = & \frac{1}{\Delta t} \sum_{r= -p + 1}^p \gamma_{p,j}^{1,0}  f(P^1_{i,j} (r \Delta t) ) + O(\Delta x^{2p})\\
& = & \frac{1}{\Delta t} \sum_{r= -p + 1}^p \gamma_{p,j}^{1,0} \sum_{k=0}^{2p-1} \frac{1}{k!}d^k (f \circ P^1_{i,j} ) (t_n) r^k \Delta t^k  + O(\Delta x^{2p-1})\\
& = & \frac{1}{\Delta t}  \sum_{k=0}^{2p-1} \frac{1}{k!}d^k (f \circ P^1_{i,j} ) (t_n)  \Delta t^k  \sum_{r= -p + 1}^p \gamma_{p,j}^{1,0} r^k + O(\Delta x^{2p-1})\\
& = & d^1 (f \circ P^1_{i,j} ) (t_n) + O(\Delta x^{2p-1})\\
& = & \partial_t f(u)(x_{i+j}, t_n) + O(\Delta x^{2p-1}),
\end{eqnarray*}
where \eqref{propbasc2} has been used. 
This result can be extended by induction to every $k $ between 1 and $2p -1$ as follows:
\begin{equation}
\tilde{f}^{(k)}_{i,j} = \partial_t^k f(u) (x_{i+j}, t_n) + O(\Delta t^{2p - k}), \quad k=1, \dots, 2p-1.
\end{equation}
Using this equality we get:
\begin{align*}
 &u(x_i, t_{n+1}) -   u(x_i, t_n) + \frac{\Delta t}{\Delta x}\left( F^p_{i+1/2} - F^p_{i-1/2} \right) \\
  &=   u(x_i, t_{n+1})  -   u(x_i, t_n) + \frac{1}{\Delta x}\sum_{k=1}^{2p}\frac{\Delta t^k}{k!} \left( A^{0,1/2}_{p,0}(\tilde{f}^{(k-1)}_{i, \bcdot}, \Delta x) - A^{0,1/2}_{p,0}(\tilde{f}^{(k-1)}_{i-1, \bcdot}, \Delta x) \right) \\
   & =  u(x_i, t_{n+1})  - u(x_i, t_n) + \frac{1}{\Delta x}\sum_{k=1}^{2p}\frac{\Delta t^k}{k!} \left( A^{0,1/2}_{p,i}(\partial_t^{k-1} f(u), \Delta x ) - A^{0,1/2}_{p,i-1}(\partial_t^{k-1} f(u), \Delta x )  \right)\\
    &\quad + O(\Delta x^{2p + 1})\\
    &= u(x_i, t_{n+1})  - u(x_i, t_n) + \frac{1}{\Delta x}\sum_{k=1}^{2p}\frac{\Delta t^k}{k!} D^{1}_{p,i}(\partial_t^{k-1} f(u), \Delta x )  + O(\Delta x^{2p + 1})\\    
     & = u(x_i, t_{n+1})  - u(x_i, t_n) + \frac{1}{\Delta x}\sum_{k=1}^{2p}\frac{\Delta t^k}{k!} \partial_t^{k-1} f(u) (x_i, t_n)  + O(\Delta x^{2p + 1})\\   
      &= u(x_i, t_{n+1})  - u(x_i, t_n) -  \frac{1}{\Delta x}\sum_{k=1}^{2p}\frac{\Delta t^k}{k!} \partial_t^{k}u (x_i, t_n)  + O(\Delta x^{2p + 1})\\   
      & = O(\Delta x^{2p + 1}). 
\end{align*}
\hfill$\Box$
\end{dem}

\textbf{Remark}:
In the Approximate Taylor method proposed in \cite{ZORIO}  the 
derivatives $\tilde u^{(k+1)}_i $ are computed by applying the $2 p_k + 
1$-point centered differentiation formula for
first derivatives to $\tilde{f}^{(k)}_{i} $, where $p_k$ is given by 
\eqref{p_k}: notice that $p_k$ decreases as $k$ increases. The same 
reduction of the stencil used to compute
$\tilde u^{(k)}_{i,j} $
could be applied here, what would allow us to reduce the number of 
computations while preserving the overall order of accuracy. 
Nevertheless, the resulting method will not be an extension of the 
linear Lax-Wendroff method. On the other hand, the CPU reduction will be 
not significant.

\subsection{Example: Fourth order compact approximate Taylor method}

Since the method is conservative, we will only show in detail  how  to compute the  numerical flux (\ref{cat1}) and, to do this, it is enough to
specify how to compute
\begin{equation} \label{flujp4}
\kappa^k_{i + 1/2} =  A^{0, 1/2}_{2,0}(\tilde{f}^{(k-1)}_{i, \bcdot} , \Delta x), \quad k =1, 2, 3, 4. 
\end{equation}

The procedure is as follows:
\begin{itemize}

\item $\kappa^{1}_{i+1/2}$: First the assignment 
  $$\tilde{f}^{(0)}_{i,j}= f(u^n_{i+j}), \quad j=-1,...,2$$ 
  is done and  then:
$$
\kappa^{1}_{i+1/2}  = A^{0, 1/2}_{2,0}(\tilde{f}^{(0)}_{i, \bcdot}, \Delta x)  
 =  \frac{ -\tilde{f}^{(0)}_{i,-1} +7 \tilde{f}^{(0)}_{i,0} +7\tilde{f}^{(0)}_{i,1}-\tilde{f}^{(0)}_{i,2}}{12}       .
$$

\item $\kappa^{2}_{i+1/2}$: The first order time derivatives of $u$ at the nodes $i-1, \dots, i + 2$ are approximated by applying the corresponding differentiation numerical formula to $\tilde{f}^{(0)}_{i,j}$:
\begin{eqnarray*}
\tilde{u}^{(1)}_{i,-1} = & -  A^{1,-1}_{2,0}(\tilde{f}^{(0)}_{i, \bcdot}, \Delta x)  = & -\frac{ 11/6 \tilde{f}^{(0)}_{i,-1} - 3 \tilde{f}^{(0)}_{i,0} +3/2 \tilde{f}^{(0)}_{i,1} - 1/3 \tilde{f}^{(0)}_{i,2}}{\Delta x},\\
\tilde{u}^{(1)}_{i,0} = & -  A^{1,0}_{2,0}(\tilde{f}^{(0)}_{i, \bcdot}, \Delta x)  = & -\frac{ 1/3 \tilde{f}^{(0)}_{i,-1} +1/2 \tilde{f}^{(0)}_{i,0} - \tilde{f}^{(0)}_{i,1} + 1/6 \tilde{f}^{(0)}_{i,2}}{\Delta x},\\
\tilde{u}^{(1)}_{i,1} = & -  A^{1,1}_{2,0}(\tilde{f}^{(0)}_{i, \bcdot}, \Delta x)  = & -\frac{-1/6 \tilde{f}^{(0)}_{i,-1} + \tilde{f}^{(0)}_{i,0} -1/2 \tilde{f}^{(0)}_{i,1} - 1/3 \tilde{f}^{(0)}_{i,2}}{\Delta x},\\
\tilde{u}^{(1)}_{i,2} = & -  A^{1,2}_{2,0}(\tilde{f}^{(0)}_{i, \bcdot}, \Delta x)  = & -\frac{ 1/3 \tilde{f}^{(0)}_{i,-1} -3/2 \tilde{f}^{(0)}_{i,0} +3 \tilde{f}^{(0)}_{i,1} -11/6 \tilde{f}^{(0)}_{i,2}}{\Delta x}.\\
\end{eqnarray*}

Next first order Taylor expansions are used to approximate the values of the flux sixteen space-time local nodes: for $r = -1, \dots, 2$
\begin{eqnarray*}
\tilde{f}^{1,r}_{i,-1} = &  f( u^{n}_{i-1} + r \Delta t \, \tilde{u}^{(1)}_{i,-1} ),\\
\tilde{f}^{1,r}_{i,0}  = & f(u^{n}_{i,0} + r \Delta t \, \tilde{u}^{(1)}_{i,0}),\\
\tilde{f}^{1,r}_{i,1} = & f( u^{n}_{i+1} + r \Delta t \, \tilde{u}^{(1)}_{i,+1}),\\
\tilde{f}^{1,r}_{i,2} = & f( u^{n}_{i+2} + r \Delta t \, \tilde{u}^{(1)}_{i,+2}).\\
\end{eqnarray*}
 Then, the first order time derivates of the flux at the nodes $i-1, \dots, i+2$  are approximated by applying the corresponding differentiation numerical formula to $\tilde{f}^{1,r}_{i,j}$:
\begin{eqnarray*}
\tilde{f}^{(1)}_{i,-1} = &  A^{1,0}_{2,n}(\tilde{f}^{1,\bcdot}_{i,-1}, \Delta t)  = & \frac{ -1/3 \tilde{f}^{1,n-1}_{i,-1} - 1/2 \tilde{f}^{1,n}_{i,-1} + \tilde{f}^{1,n+1}_{i,-1} - 1/6 \tilde{f}^{1,n+2}_{i,-1}}{\Delta t},\\
\tilde{f}^{(1)}_{i,0} = &  A^{1,0}_{2,n}(\tilde{f}^{1,\bcdot}_{i,0}, \Delta t)  = & \frac{ -1/3 \tilde{f}^{1,n-1}_{i,0} - 1/2 \tilde{f}^{1,n}_{i,0} + \tilde{f}^{1,n+1}_{i,0} - 1/6 \tilde{f}^{1,n+2}_{i,0}}{\Delta t},\\
\tilde{f}^{(1)}_{i,1} = &  A^{1,0}_{2,n}(\tilde{f}^{1,\bcdot}_{i,1}, \Delta t)  = & \frac{ -1/3 \tilde{f}^{1,n-1}_{i,1} - 1/2 \tilde{f}^{1,n}_{i,1} + \tilde{f}^{1,n+1}_{i,1} - 1/6 \tilde{f}^{1,n+2}_{i,1}}{\Delta t},\\
\tilde{f}^{(1)}_{i,2} = &  A^{1,0}_{2,n}(\tilde{f}^{1,\bcdot}_{i,2}, \Delta t)  = & \frac{ -1/3 \tilde{f}^{1,n-1}_{i,2} - 1/2 \tilde{f}^{1,n}_{i,2} + \tilde{f}^{1,n+1}_{i,2} - 1/6 \tilde{f}^{1,n+2}_{i,2}}{\Delta t}.\\
\end{eqnarray*}
Finally;
$$
\kappa^{2}_{i+1/2}  = A^{0, 1/2}_{2,0}(\tilde{f}^{(1)}_{i, \bcdot}, \Delta x)  = \frac{ -\tilde{f}^{(1)}_{i,-1} +7 \tilde{f}^{(1)}_{i,0} +7\tilde{f}^{(1)}_{i,1}-\tilde{f}^{(1)}_{i,2}}{12}       .
$$
 
\item $\kappa^{3}_{i+1/2}$: the second order time derivatives at the nodes are approximated by
\begin{eqnarray*}
\tilde{u}^{(2)}_{i,-1} = & -  A^{1,-1}_{2,0}(\tilde{f}^{(1)}_{i, \bcdot}, \Delta x)  = & -\frac{ 11/6 \tilde{f}^{(1)}_{i,-1} - 3 \tilde{f}^{(1)}_{i,0} +3/2 \tilde{f}^{(1)}_{i,1} - 1/3 \tilde{f}^{(1)}_{i,2}}{\Delta x},\\
\tilde{u}^{(2)}_{i,0} = & -  A^{1,0}_{2,0}(\tilde{f}^{(1)}_{i, \bcdot}, \Delta x)  = & -\frac{ 1/3 \tilde{f}^{(1)}_{i,-1} +1/2 \tilde{f}^{(1)}_{i,0} - \tilde{f}^{(1)}_{i,1} + 1/6 \tilde{f}^{(1)}_{i,2}}{\Delta x},\\
\tilde{u}^{(2)}_{i,1} = & -  A^{1,1}_{2,0}(\tilde{f}^{(1)}_{i,\bcdot}, \Delta x)  = & -\frac{-1/6 \tilde{f}^{(1)}_{i,-1} + \tilde{f}^{(1)}_{i,0} -1/2 \tilde{f}^{(1)}_{i,1} - 1/3 \tilde{f}^{(1)}_{i,2}}{\Delta x},\\
\tilde{u}^{(2)}_{i,2} = & -  A^{1,2}_{2,0}(\tilde{f}^{(1)}_{i, \bcdot}, \Delta x)  = & -\frac{ 1/3 \tilde{f}^{(1)}_{i,-1} -3/2 \tilde{f}^{(1)}_{i,0} +3 \tilde{f}^{(1)}_{i,1} -11/6 \tilde{f}^{(1)}_{i,2}}{\Delta x}.\\
\end{eqnarray*}
Second order Taylor expansions are used to compute the fluxes at the sixteen nodes in the space-time mesh: for $r = -1, \dots 2$
\begin{eqnarray*}
\tilde{f}^{2,r}_{i,-1} = &  f\left(u^{n}_{i-1} + r \Delta t \, \tilde{u}^{(1)}_{i,-1}+ \frac{r^2 \Delta t^2}{2} \, \tilde{u}^{(2)}_{i,-1}\right),\\
\tilde{f}^{2,r}_{i,0}  = & f\left( u^{n}_{i} + r \Delta t \, \tilde{u}^{(1)}_{i,0}+ \frac{r^2 \Delta t^2}{2} \, \tilde{u}^{(2)}_{i,0}]\right),\\
\tilde{f}^{2,r}_{i,1} = & f\left(u^{n}_{i+1} + r \Delta t \, \tilde{u}^{(1)}_{i,1}+ \frac{r^2 \Delta t^2}{2} \, \tilde{u}^{(2)}_{i,1}\right),\\
\tilde{f}^{2,r}_{i,2} = & f\left( u^{n}_{i+2} + r \Delta t \, \tilde{u}^{(1)}_{i,2}+ \frac{r^2 \Delta t^2}{2} \, \tilde{u}^{(2)}_{i,2}\right).\\
\end{eqnarray*}
Next, compute
\begin{eqnarray*}
\tilde{f}^{(2)}_{i,-1} = &  A^{2,0}_{2,n}(\tilde{f}^{2,\bcdot}_{i,-1}, \Delta t)  = & \frac{ \tilde{f}^{2,n-1}_{i,-1} - 2 \tilde{f}^{2,n}_{i,-1} + \tilde{f}^{2,n+1}_{i,-1} }{\Delta t^2},\\
\tilde{f}^{(2)}_{i,0} = &  A^{2,0}_{2,n}(\tilde{f}^{2,\bcdot}_{i, 0}, \Delta t)  = & \frac{ \tilde{f}^{2,n-1}_{i,0} -2 \tilde{f}^{2,n}_{i,0} + \tilde{f}^{2,n+1}_{i,0} }{\Delta t^2},\\
\tilde{f}^{(2)}_{i,1} = &  A^{2,0}_{2,n}(\tilde{f}^{2,\bcdot}_{i,1}, \Delta t)  = & \frac{  \tilde{f}^{2,n-1}_{i,1} - 2 \tilde{f}^{2,n}_{i,1} + \tilde{f}^{2,n+1}_{i,1} }{\Delta t^2},\\
\tilde{f}^{(2)}_{i,2} = &  A^{2,0}_{2,n}(\tilde{f}^{2,\bcdot}_{i,2}, \Delta t)  = & \frac{  \tilde{f}^{2,n-1}_{i,2} - 2 \tilde{f}^{2,n}_{i,2} + \tilde{f}^{2,n+1}_{i,2} }{\Delta t^2}.\\
\end{eqnarray*}
And finally;
$$
\kappa^{3}_{i+1/2} =  A^{0, 1/2}_{2,0}(\tilde{f}^{(2)}_{i, \bcdot}, \Delta x)   =  \frac{ -\tilde{f}^{(2)}_{i,-1} +7 \tilde{f}^{(2)}_{i,0} +7\tilde{f}^{(2)}_{i,1}-\tilde{f}^{(2)}_{i,2}}{12}       .
$$
 
\item$\kappa^{4}_{i+1/2}$: the third order time derivatives at the nodes are approximated by
\begin{eqnarray*}
\tilde{u}^{(3)}_{i,-1} = & -  A^{1,-1}_{2,0}(\tilde{f}^{(2)}_{i, \bcdot}, \Delta x)  = &- \frac{ 11/6 \tilde{f}^{(2)}_{i,-1} - 3 \tilde{f}^{(2)}_{i,0} +3/2 \tilde{f}^{(2)}_{i,1} - 1/3 \tilde{f}^{(2)}_{i,2}}{\Delta x},\\
\tilde{u}^{(3)}_{i,0} = & -  A^{1,0}_{2,0}(\tilde{f}^{(2)}_{i, \bcdot}, \Delta x)  = &- \frac{ 1/3 \tilde{f}^{(2)}_{i,-1} +1/2 \tilde{f}^{(2)}_{i,0} - \tilde{f}^{(2)}_{i,1} + 1/6 \tilde{f}^{(2)}_{i,2}}{\Delta x},\\
\tilde{u}^{(3)}_{i,1} = & -  A^{1,1}_{2,0}(\tilde{f}^{(2)}_{i, \bcdot}, \Delta x)  = &- \frac{-1/6 \tilde{f}^{(2)}_{i,-1} + \tilde{f}^{(2)}_{i,0} -1/2 \tilde{f}^{(2)}_{i,1} - 1/3 \tilde{f}^{(2)}_{i,2}}{\Delta x},\\
\tilde{u}^{(3)}_{i,2} = & -  A^{1,2}_{2,0}(\tilde{f}^{(2)}_{i, \bcdot}, \Delta x)  = & -\frac{ 1/3 \tilde{f}^{(2)}_{i,-1} -3/2 \tilde{f}^{(2)}_{i,0} +3 \tilde{f}^{(2)}_{i,1} -11/6 \tilde{f}^{(2)}_{i,2}}{\Delta x}.\\
\end{eqnarray*}
Compute the approximations of the fluxes: for $r = -1, \dots, 2$
\begin{eqnarray*}
\tilde{f}^{3,r}_{i,-1} = &  f\left( u^{n}_{i-1} + r \Delta t \, \tilde{u}^{(1)}_{i,-1}+ \frac{r^2 \Delta t^2}{2} \, \tilde{u}^{(2)}_{i,-1}+ \frac{r^3 \Delta t^3}{6} \, \tilde{u}^{(3)}_{i,-1}\right),\\
\tilde{f}^{3,r}_{i,0}  = & f\left( u^{n}_{i} + r \Delta t \, \tilde{u}^{(1)}_{i,0}+ \frac{r^2 \Delta t^2}{2} \, \tilde{u}^{(2)}_{i,0}+ \frac{r^3 \Delta t^3}{6} \, \tilde{u}^{(3)}_{i,0}\right),\\
\tilde{f}^{3,r}_{i,1} = & f\left(u^{n}_{i+1} + r \Delta t \, \tilde{u}^{(1)}_{i,1}+ \frac{r^2 \Delta t^2}{2} \, \tilde{u}^{(2)}_{i,1}+\frac{r^3 \Delta t^3}{6} \, \tilde{u}^{(3)}_{i,1}\right),\\
\tilde{f}^{3,r}_{i,2} = & f\left( u^{n}_{i+2} + r \Delta t \, \tilde{u}^{(1)}_{i,2}+ \frac{r^2 \Delta t^2}{2} \, \tilde{u}^{(2)}_{i,2}+\frac{r^3 \Delta t^3}{6} \, \tilde{u}^{(3)}_{i,2}\right).\\
\end{eqnarray*}
Next, compute:
\begin{eqnarray*}
\tilde{f}^{(3)}_{i,-1} = &  A^{3,0}_{2,n}(\tilde{f}^{3,\bcdot}_{i,-1}, \Delta t)  = & \frac{ -\tilde{f}^{3,n-1}_{i,-1} +3 \tilde{f}^{3,n}_{i,-1} -3 \tilde{f}^{3,n+1}_{i,-1} +\tilde{f}^{3,n+2}_{i,-1} }{\Delta t^3},\\
\tilde{f}^{(3)}_{i,0} = &  A^{3,0}_{2,n}(\tilde{f}^{3,\bcdot}_{i,0}, \Delta t)  = & \frac{ -\tilde{f}^{1,n-1}_{i,0} +3 \tilde{f}^{3,n}_{i,0} -3\tilde{f}^{3,n+1}_{i,0} +\tilde{f}^{3,n+2}_{i,0}}{\Delta t^3},\\
\tilde{f}^{(3)}_{i,1} = &  A^{3,0}_{2,n}(\tilde{f}^{3,\bcdot}_{i,1}, \Delta t)  = & \frac{  -\tilde{f}^{3,n-1}_{i,1} +3 \tilde{f}^{3,n}_{i,1} -3 \tilde{f}^{3,n+1}_{i,1}+ \tilde{f}^{3,n+2}_{i,1} }{\Delta t^3},\\
\tilde{f}^{(3)}_{i,2} = &  A^{3,0}_{2,n}(\tilde{f}^{3,\bcdot}_{i,2}, \Delta t)  = & \frac{  -\tilde{f}^{3,n-1}_{i,2} +3 \tilde{f}^{3,n}_{i,2} -3 \tilde{f}^{3,n+1}_{i,2}+\tilde{f}^{3,n+2}_{i,2} }{\Delta t^3}.\\
\end{eqnarray*}
Finally;
$$
\kappa^{4}_{i+1/2}  = A^{0, 1/2}_{2,0}(\tilde{f}^{(2)}_{i, \bcdot}, \Delta x) = \frac{ -\tilde{f}^{(2)}_{i,-1} +7 \tilde{f}^{(2)}_{i,0} +7\tilde{f}^{(2)}_{i,1}-\tilde{f}^{(2)}_{i,2}}{12}       .
$$

\end{itemize}

If  $f(u)=au$,  then: 
\begin{eqnarray}\label{flujo_cuarto_orden} 
F^2_{i+1/2} & = & \frac{a }{12} ( -u^n_{i-1}+7 u^n_{i} +7 u^n_{i+1} -u^n_{i+2}) +  \frac{a^2 \Delta t}{24 \Delta x}( -u^n_{i-1}+15 u^n_{i} -15 u^n_{i+1} +u^n_{i+2}) \\
& & + \frac{a^3 \Delta t^2}{12 \Delta x^2}( u^n_{i-1}- u^n_{i} - u^n_{i+1} + u^n_{i+2}) + \frac{a^4 \Delta t^3}{24 \Delta x^3} ( u^n_{i-1}-3 u^n_{i} +3 u^n_{i+1}-u^n_{i+2}),  
\end{eqnarray}
which coincides with the numerical flux of the  fourth order Lax-Wendroff in conservative form.

\section{Shock-capturing techniques}
 Although the Compact Approximate Taylor methods are linearly  stable in the $L^2$ sense under the usual $CFL$-1 condition, they may  produce strong oscillations close to a discontinuity of the solution. The goal of this section is to modify the numerical method to avoid these oscillations. Two different techniques are considered here:

\subsection{FLUX LIMITER-CAT methods}
 We consider the numerical method \eqref{cons} with
\begin{equation}\label{centralflux}
F_{i+1/2}= (1 - \varphi_{i+1/2} ) F^{L}_{i+1/2} + \varphi_{i+1/2} F^{p}_{i+1/2},
\end{equation}
where $ F^{L}_{i+1/2}$ is a first order robust numerical flux,\, $F^{p}_{i+1/2}$ is given by \eqref{cat1}, and $ \varphi_{i+1/2}$ is the flux limiter function, see \cite{kemm} \cite{book_leveque_2002}, \cite{toro2009riemann}. We consider here 
\begin{equation}\label{indicador1}
\varphi_{i+1/2}  =  \varphi(r_{i+1/2}) , 
\end{equation}
where $\varphi$ is the van Albada second version flux limiter:
\begin{equation}\label{indicator}
\varphi(r)= \max \left(0,\frac{2 r} { 1+r^2} \right),
\end{equation}
and
$$
r_{i+1/2}=\frac{\Delta upw}{\Delta loc}
= \left\{
\begin{array}{cl}
\displaystyle \frac{u^n_i-u^n_{i-1}}{u^n_{i+1}-u^n_{i}} & \mbox {if } a_{i+1/2}  >0, \\
\displaystyle \frac{u^n_{i+2}-u^n_{i+1}}{u^n_{i+1}-u^n_{i}} & \mbox {if } a_{i+1/2}  < 0,
\end{array}\right.
$$
where, $a_{i+1/2}$ is an estimate of the wave speed.

\subsection{WENO-CAT methods}

Following \cite{JIANXIAN1064827502412504} WENO reconstructions of the flux are used to stabilize the method. The only differences with the algorithm described in Section \ref{sec:CAT}
are the computation of $ \tilde{u}  ^{(1)}_{i,j} $, that is now  performed as follows:
$$
\tilde{u}  ^{(1)}_{i,j} = - \frac{\hat f_{i + j +  1/2} - \hat f_{i + j -1/2}}{\Delta x},
$$
where $\hat f_{i + 1/2}$ denotes the WENO flux splitting, reconstructions at $x_{i + 1/2}$  of the flux function described in \cite{Shu1989}. The expression of the numerical flux is then given by:
\begin{equation}\label{WENOCAT}
F^p_{i+1/2}  = \hat f_{i + 1/2} + \sum_{k=2}^{2p} \frac{\Delta t^{k-1}}{k!}A^{0, 1/2}_{p,0}(\tilde{f}_{i, \bcdot}^{(k-1)}, \Delta x).
\end{equation}

\subsection{Systems of conservation laws}
Although the methods have been described in the scalar case for simplicity, they can be easily applied to systems using vector notation. 

\section{Numerical Experiments}
In this section we apply to some  scalar conservation laws and to the  1D Euler equation  the following numerical methods:
\begin{itemize}
\item LW-CAT$q$: Compact Approximate Taylor method of order $q$ (space and time).
\item FL-CAT$q$: Compact Approximate Taylor method of order $q$ with flux limiter technique. The first order methods  considered are Lax-Friedrich for scalar problems and HLL for systems. 
\item WENO$s$-CAT$q$:  Compact Approximate Taylor method of order $q$ with WENO reconstructions of order $s$ to compute $\tilde{u}  ^{(1)}_{t,i}$.
\item WENO$s$-RK$q$: WENO method of order $s$ for the space discretization and TVD-RK$q$ for the time discretization, see \cite{Shu1989}.
\item WENO$s$-LWA$q$:  Approximate Taylor method of order $q$ with WENO reconstructions of order $s$ to compute $\tilde{u}  ^{(1)}_{t,i}$, see \cite{ZORIO}.
\end{itemize}

\subsection{Linear transport equation}
We consider first \eqref{edplinsc} with $a = 1$, in the space interval  $[0,1]$, with initial condition
\begin{equation} \label{square_step_test} 
u(x,0)= \left\{ \begin{array}{l l} 1 & 0 \leq x< 1/2, \\
2 & 1/2 \leq x < 1,
\end{array}\right.
\end{equation}
  and periodic boundary conditions. A uniform mesh with $N = 80$ points is considered and the LW-CAT method (that, in this case, coincides with the Lax-Wendroff method) is applied for  $p=1,\dots,5$.
  
\begin{figure}[!ht]
	\small
	\setlength{\unitlength}{1mm}
	\centering	
	\begin{picture}(70,75)
	\put(-17,3){\makebox(100,65)[c]{
	\includegraphics[width=18cm,height=8cm]{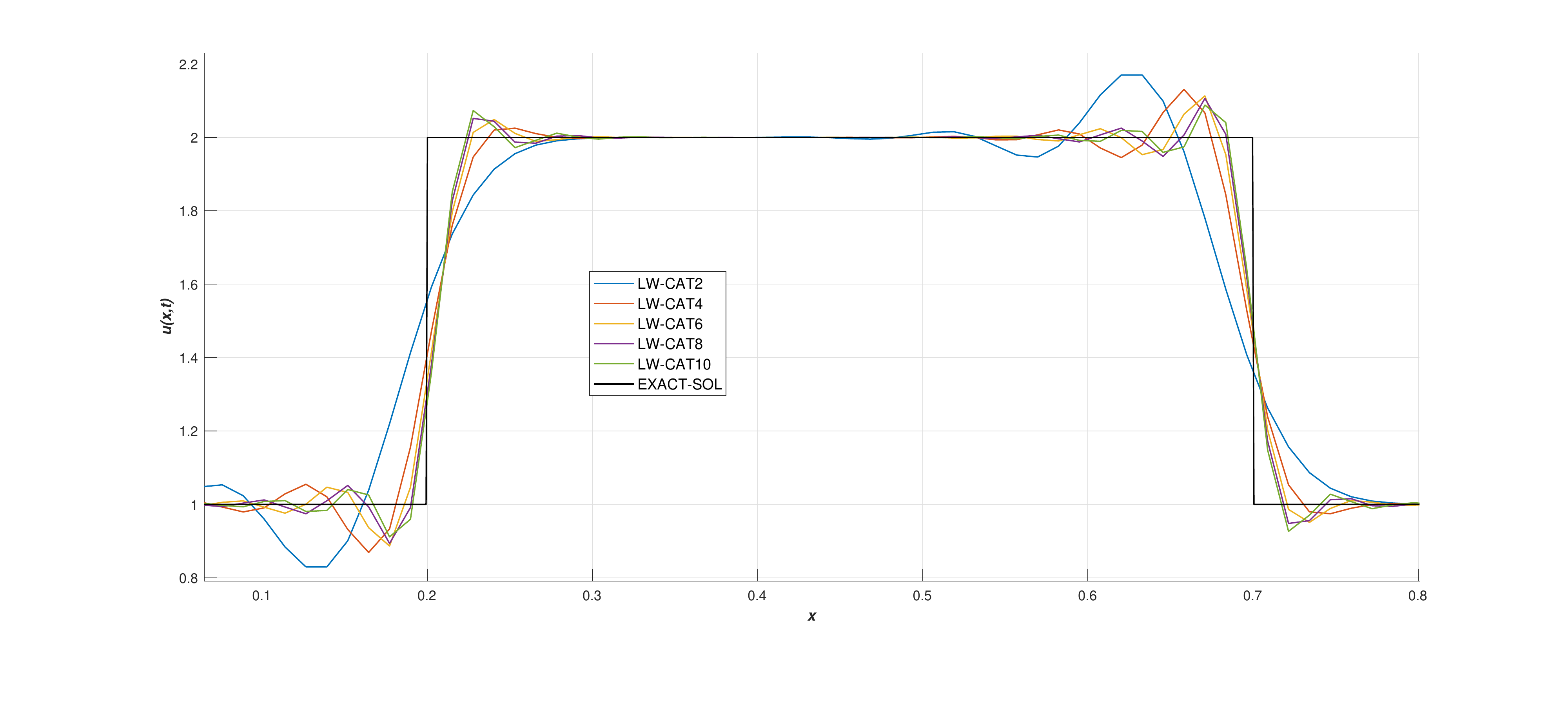}}}
	\end{picture}
	\vspace{-1.2 cm}
	\caption{Transport equation with initial condition  \eqref{square_step_test}, $CFL=0.9$ and $t=1$s: numerical results obtained with LW-CATq, $q = 2, 4, 6, 8, 10$.}
	\label{Num_02}
	\end{figure}

\begin{figure}[!ht]
	\small
	\setlength{\unitlength}{1mm}
	\centering	
	\begin{picture}(70,85)
	\put(-17,3){\makebox(100,65)[c]{
	\includegraphics[width=18cm,height=10cm]{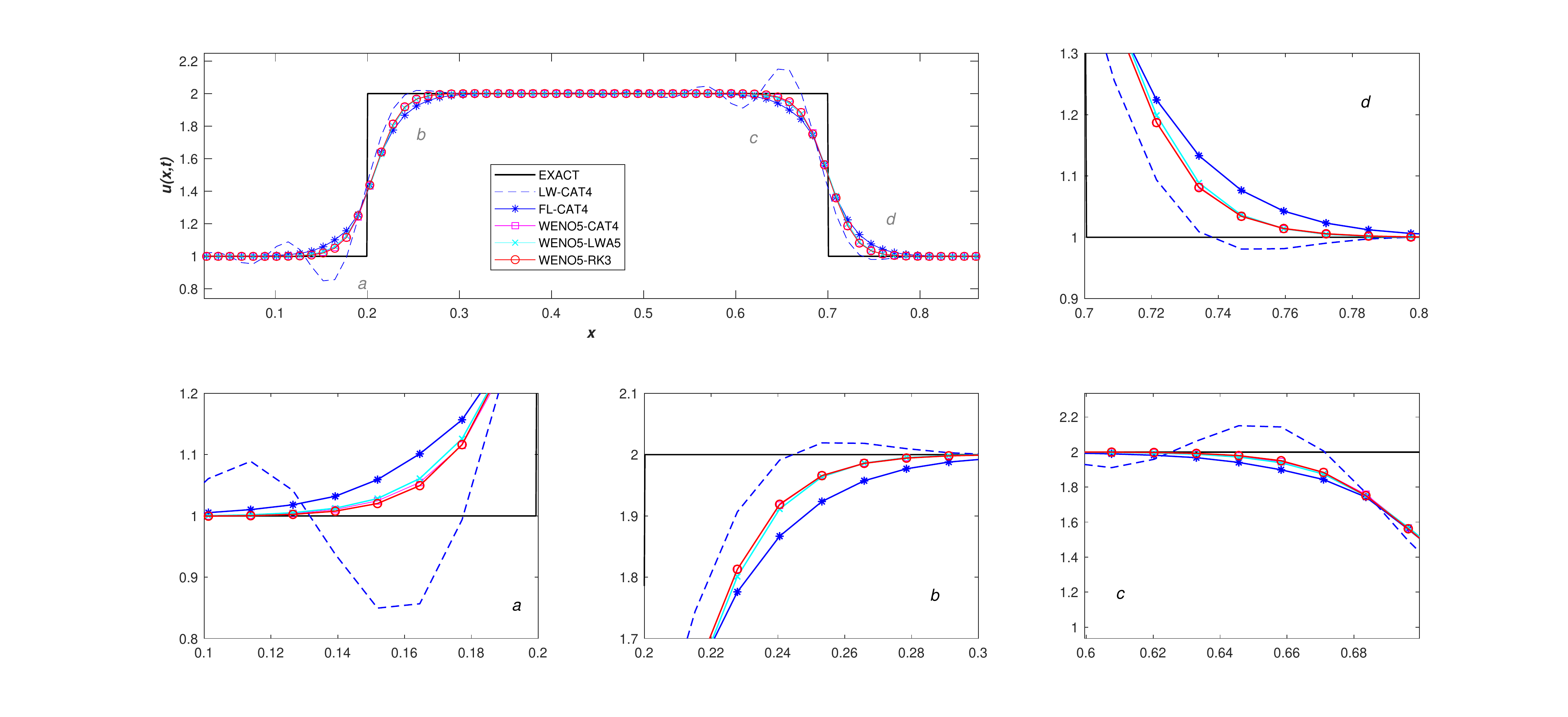}}}
	\end{picture}
	\vspace{0.5 cm}
	\caption{Transport equation with initial condition (\ref{square_step_test}),  $CFL=0.5$ and $t=1$s.  Left-up: general viewt. \textbf{a},\textbf{b},\textbf{c} and \textbf{d}: enlarged view of interest areas. }
	\label{Num_03}
	\end{figure}

Numerical simulations  are shown in Figure \ref{Num_02}: the $L^2$ stability of the scheme and the appearance of oscillations near the discontinuities can be observed. Next, we apply to the same problem  LW-CAT4, FL-CAT4, WENO5-CAT4, WENO5-RK3, and WENO5-LWA5 methods. A general view is shown in Figure \ref{Num_03} together with a zoom of the area of interest. 
As it can be observed, the results  given by WENO5-CAT4, WENO5-RK3 and WENO5-LWA5 are almost identical. Nevertheless, as it will be seen in the next test problem, WENO5-CAT4 still gives good results for CFL close to one, what is not the case for WENO5-RK3 or WENO5-LWA5.

Using the algorithm in Section 3.2 CAT methods are easily extended to any $2p$ order, nevertheless, this involves a significant increase of the CPU time simulation and flops (number of operations required), that should be considered, see Table \ref{Tabla_exp02}.

\begin{table}[htbp]
\begin{center}
\begin{tabular}{|c| c c c c c c c c c|}
\hline
\backslashbox{Rate}{Order}  & LW-CAT2  &   &LW-CAT4 &   &LW-CAT6 &   &LW-CAT8&   &LW-CAT10    \\
\hline 
 time  & 1 &   &2.98&   &7.72&   &18.87&   &42.66    \\
\hline 
 flops & 1 &   &1.61&   &2.51&   &3.69&   &5.16    \\
\hline
\end{tabular}
\caption{Average rate time and flops to increase from  LW-CAT2 to LW-CAT2p, with $p=2,3,4,5$ for the scalar transport equation  with initial conditions (\ref{square_step_test}) and  $t=1$s.}
\label{Tabla_exp02}
\end{center}
\end{table}

Finally, we consider \eqref{edplinsc} in the space interval $[0,2]$ with initial condition, 
\begin{equation} \label{sin_test}
u(x,0)=0.25 \sin(\pi x),
\end{equation}
and periodic boundary conditions. Table \ref{Tabla_exp03} shows the error and the empirical order for  LW-CAT2, LW-CAT4, LW-CAT6, and Table \ref{Tabla_exp04} \,  for WENO5-RK3 and WENO5-LWA5 which coincides in all cases with the theoretical one.  For smooth solutions WENO-CATq and FL-CATq reduce to the corresponding LW-CATq, so that the accuracy test is not necessary.

\begin{table}[htbp]
\begin{center}
\begin{tabular}{|l| l l| l l| l l|}
\hline
  &LW-CAT2         &      & LW-CAT4       &        &LW-CAT6   &    \\
\hline
  $\Delta x$    &Error $\Vert \cdot \Vert_1 $    &Order $\Vert \cdot \Vert_1 $    &Error $\Vert \cdot \Vert_1 $    &Order $\Vert \cdot \Vert_1 $     &Error $\Vert \cdot \Vert_1 $    &Order $\Vert \cdot \Vert_1 $  \\
\hline \hline
 0.1053         &3.68e-02    &       &1.40e-02   &    	&7.88e-03    & \\
 0.0526         &6.84e-03    &2.43   &3.50e-05   &8.64 	&4.25e-08    &7.50\\
 0.0263         &1.70e-03    &2.00   &2.19e-06   &4.00 	&6.49e-10    &6.03\\
 0.0132       &4.27e-04    &2.00   &1.36e-07   &4.00 	&9.89e-12    &6.04\\
 0.0066        &1.06e-04    &2.00   &8.55e-09   &4.00 	&1.53e-13    &6.01\\
 0.0033        &2.66e-05    &2.00   &5.34e-10   &4.00 	&2.64e-15    &5.96\\ \hline
\end{tabular}
\caption{Linear transport equation with initial condition \eqref{sin_test}, $CFL=0.5$ and $t=1$s: $L^1$ errors and accuracy order for LW-CAT2p, $p=1, 2, 3$.}
\label{Tabla_exp03}
\end{center}
\end{table}

\begin{table}[htbp]
\begin{center}
\begin{tabular}{|l| l l| l l|}
\hline
                 &WENO5-RK3         &      &WENO-LWA5    &   \\
\hline
  $\Delta x$      &Error $\Vert \cdot \Vert_1 $    &Order $\Vert \cdot \Vert_1 $    &Error $\Vert \cdot \Vert_1 $    &Order $\Vert \cdot \Vert_1 $   \\
\hline \hline
 0.1053     &2.03e-03    &       &5.44e-05   &    \\
 0.0526     &6.06e-05    &5.06   &1.65e-06   &5.04  \\
 0.0263     &1.87e-06    &5.02   &5.04e-08   &5.04 \\
 0.0132     &5.83e-08    &5.00   &1.51e-09   &5.05   \\
 0.0066     &1.82e-09    &5.00   &4.41e-11   &5.10  \\
 0.0033     &5.65e-11    &5.01   &1.15e-12   &5.25  \\ \hline
\end{tabular}
\caption{Linear transport equation with initial condition \eqref{sin_test}, $CFL=0.5$ and $t=1$s: $L^1$ errors and accuracy order for WENO5-RK3 and WENO5-LWA5.}
\label{Tabla_exp04}
\end{center}
\end{table}

\subsection{Burgers equation}

\begin{figure}[!ht]
	\small
	\setlength{\unitlength}{1mm}
	\centering	
	\begin{picture}(70,65)
	\put(-20,3){\makebox(100,50)[c]{
	\includegraphics[width=19cm,height=8cm]{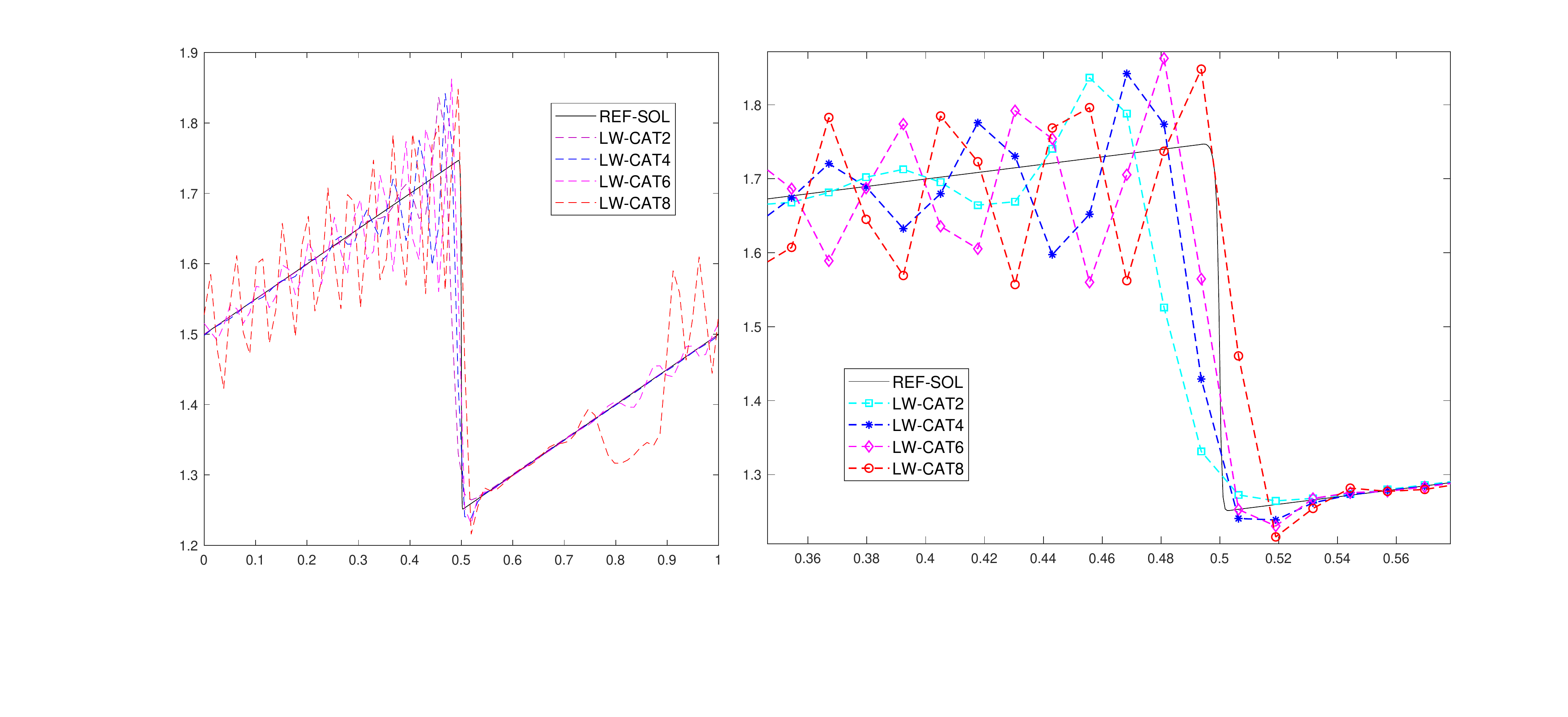}}}
	\end{picture}
	\vspace{-0.5 cm}
	\caption{Burgers equation with initial condition  \eqref{square_step_test}, $CFL=0.5$ and $t=2$s: numerical results obtained with LW-CATq, $q = 2, 4, 6, 8$. Left: 
general view. Right: enlarged view.}
	\label{Num_04}
	\end{figure}

We consider Burgers equation, i.e. \eqref{scl} with
$$
f(u) = \frac{u^2}{2}.
$$
When CAT methods are applied to approximate a discontinuous solution of this nonlinear problem, the oscillations appearing close to the shocks tend to grow and to spoil the numerical solution. 
Nevertheless, it is still possible to apply these methods by reducing the $CFL$ parameter (the reduction increases with $p$): for instance, Figure \ref{Num_04} shows the results obtained with CAT-LW2$p$,  $p=1, 2, 3 ,4$ and   $CFL={0.8,0.4,0.2,0.1} $, respectively, with initial conditions \eqref{square_step_test}. 

\begin{figure}[!ht]
	\small
	\setlength{\unitlength}{1mm}
	\centering	
	\begin{picture}(70,65)
	\put(-17,3){\makebox(100,50)[c]{
	\includegraphics[width=15cm,height=8cm]{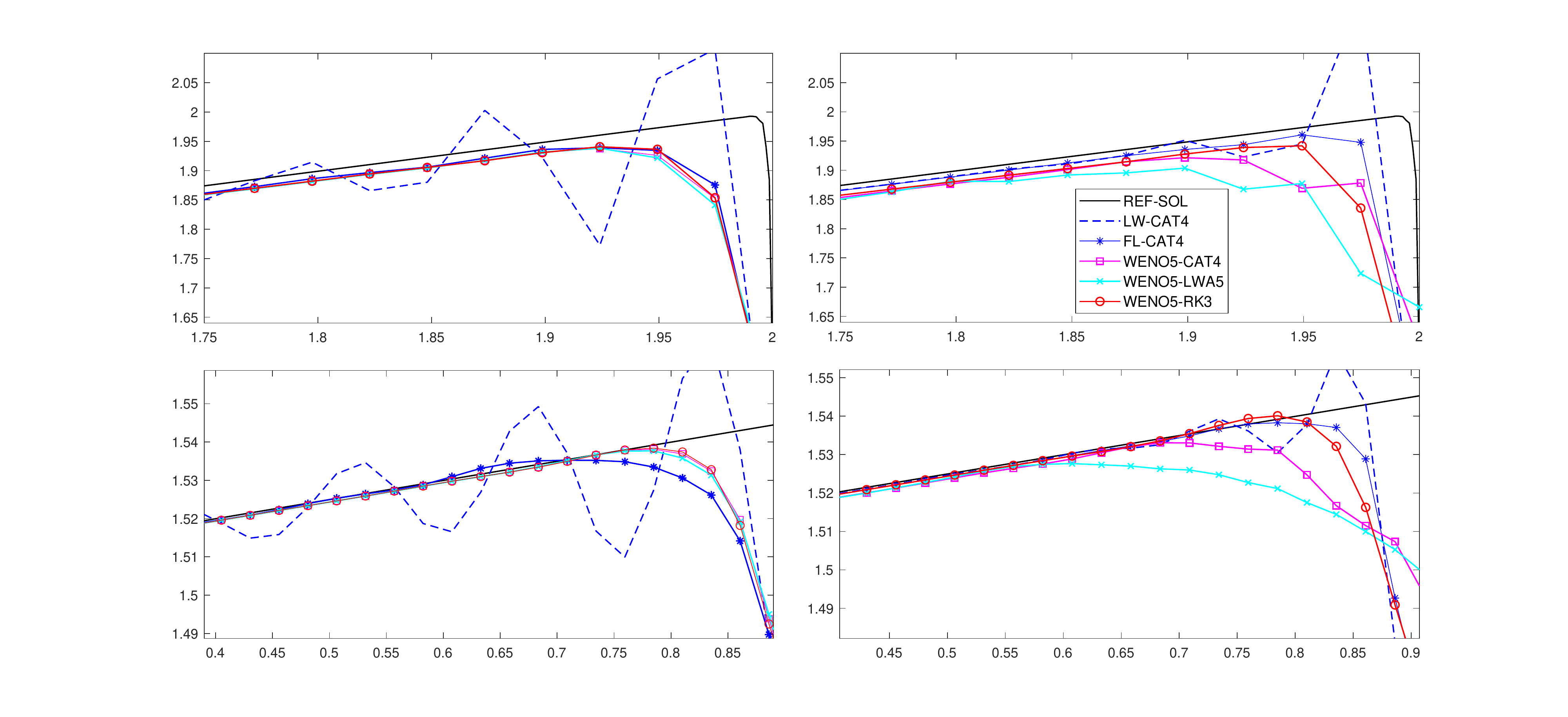}}}
	\end{picture}
	\vspace{0.5 cm}
	\caption{ Burgers equation with initial condition (\ref{square_step_test}),  $CFL=0.5$ and $CFL=0.9$, $t=2$s and $t=20$s: zooms of the numerical results.  Left-up: $CFL=0.5$ and $t=2$s. Left-down: $CFL=0.5$ and $t=20$s.  Right-up: $CFL=0.9$ and $t=2$s. Right-down: $CFL=0.9$ and $t=20$s. }
	\label{Num_06}
	\end{figure}

Next, the same test problem is solved using LW-CAT4, FL-CAT4, WENO5-RK3 and WENO5-LWA5  methods. Using $CFL = 0.5$ we obtain  numerical solutions without spurious oscillations for all the methods. Figure \ref{Num_05} shows a general view  of solutions and the van Albada flux limiter function on every inter cell. In Figure \ref{Num_06} the results are compared with those obtained with $CFl=0.9$. From Figures \ref{Num_05}, \ref{Num_06}  we can conclude:

\begin{figure}[!ht]
	\small
	\setlength{\unitlength}{1mm}
	\centering	
	\begin{picture}(70,65)
	\put(-17,3){\makebox(100,50)[c]{
	\includegraphics[width=15cm,height=8cm]{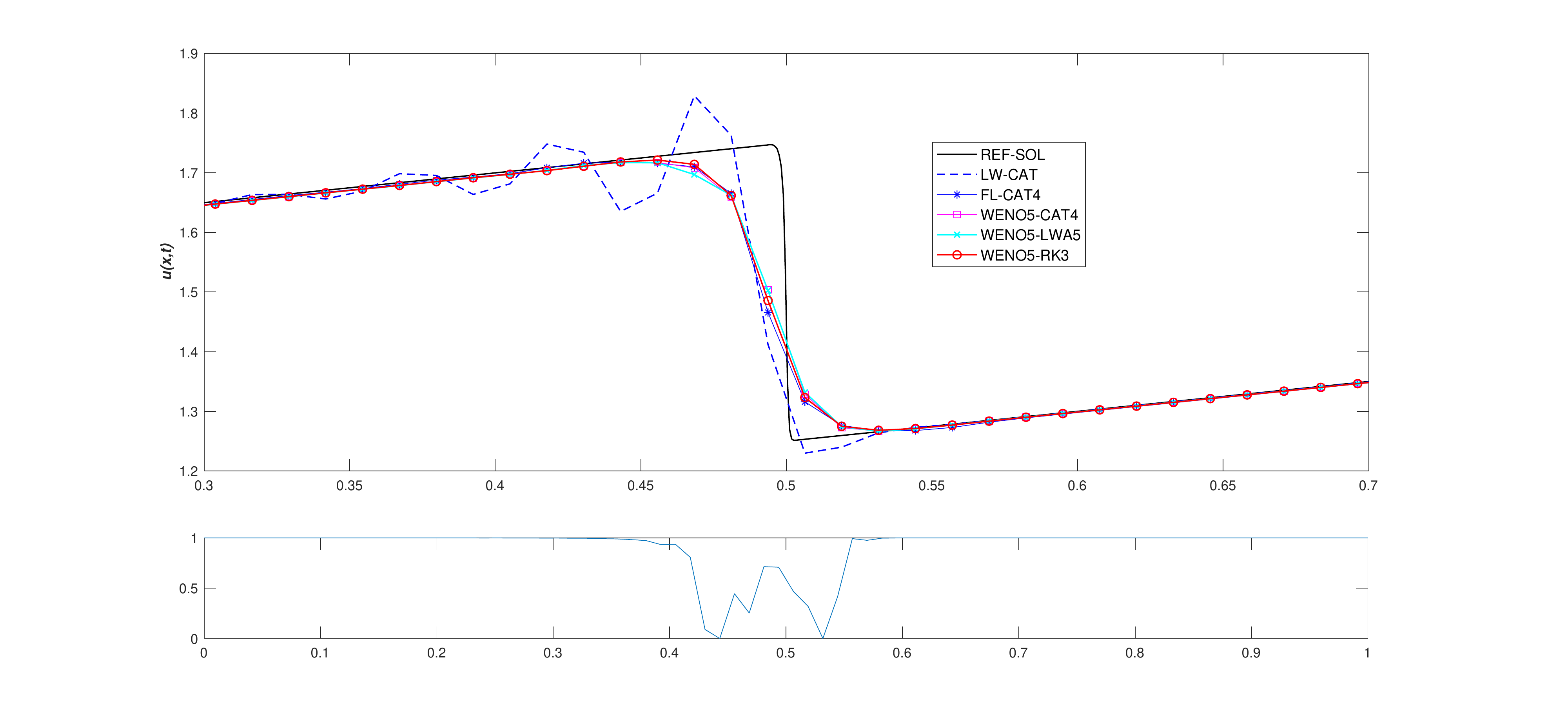}}}
	\end{picture}
	\vspace{0.5 cm}
	\caption{Burgers equation with initial condition  \eqref{square_step_test}, $CFL=0.5$ and $t=2s$.  Up: general view. Down: flux limiter function $\varphi_{i+1/2}$ for FL-CAT4.}
	\label{Num_05}
	\end{figure}

\begin{itemize}
\item  $CFL \leq 0.5$

	\begin{itemize}
     \item LW-CAT4  show oscillations near the discontinuities, but it is stable.    
	 \item FL-CAT4   is very diffusive near to the discontinuities, due to the selected  first order accurate flux limiter function.    
	\item WENO5-CAT4, WENO5-LWA4  and WENO5-RK3 show good results, stable and essentially the same values.
	\end{itemize}

\item $CFL > 0.5$

	\begin{itemize}
     \item LW-CAT4:  the amplitude of oscillations  increases near the discontinuities. However, they remain stable.    
	 \item FL-CAT4:  conversely to the previous $CFL$ condition, it shows acceptable solutions near the discontinuities.    
	\item WENO5-CAT4 ,WENO5-LWA5  and WENO5-RK3 : slight oscillations appear near the discontinuities at the beginning of the simulations. Nevertheless, as the time increases, these oscillations tend to dismiss and the result remains acceptable and stable  for WENO5-CAT4  while the solutions given by  WENO5-LWA5 is very diffusive and the one given by WENO5-RK3 is  overdamped.
	\end{itemize}
\end{itemize} 
Although FL-CAT4 shows better results for bigger $CFL$, it fails in smooth regions close to critical points and for systems (as it will be seen in Euler 1D equation). For $CFL \leq 1/2$, WENO5-LWA5 is  faster than CAT methods (in the computational cost sense): a  comparsion is shown in Table \ref{Tabla_exp_05}. However,  we can obtain good solutions for CAT methods using $CFL > 0.5$ and obtain similar CPU time computation, see Table \ref{Tabla_exp_06}.    

\begin{table}[htbp]
\begin{center}
\begin{tabular}{|c| c c c c c  c|}
\hline
Nodes  & LW-CAT4 	& FL-CAT4   &   & WENO5-CAT4 &   &WENO5-LWA5    \\
\hline 
 80    & 0.469 	& 0.4702    &   &0.5508      &   &0.1998    \\
\hline 
 160   & 0.8090    	& 0.9668    &   &0.9700      &   &0.4616     \\
\hline
 320   & 1.9838     & 1.9908    &   &2.0064      &   &1.0346     \\
\hline
\end{tabular}
\caption{ Burgers equation with initial condition (\ref{square_step_test}), $CFL=0.5$ and $t=2$s: computational time  in seconds}
\label{Tabla_exp_05}
\end{center}
\end{table}

\begin{table}[htbp]
\begin{center}
\begin{tabular}{|c| c c c|}
\hline
Nodes  & FL-CAT4  &   & WENO5-CAT4     \\
\hline 
 80    & 0.2003    &   &0.2078        \\
\hline 
 160   & 0.4024    &   &0.5001        \\
\hline
 320   & 1.0494    &   & 1.045         \\
\hline
\end{tabular}
\caption{Burgers equation with initial condition (\ref{square_step_test}), $CFL=0.9$ and $t=2$s: computational time  in seconds}
\label{Tabla_exp_06}
\end{center}
\end{table}

In order to study the order of convergence, we consider again initial condition \eqref{sin_test} and periodic boundary conditions. A reference solution at time $t =0.5$s (when the solution is still smooth) is obtained with WENO5-RK3 using  a fine grid of $1 400$ nodes.
The errors and the empirical order are shown in Table \ref{Tabla_exp07}: the numerical results verify the theoretical analysis. 
  
\begin{table}[htbp]
\begin{center}
\begin{tabular}{|l| l l| l l| l l|}
\hline
                  &LW-CAT2         &      & LW-CAT4       &        &LW-CAT6   &    \\
\hline
  $\Delta x$  &Error $\Vert \cdot \Vert_1 $    &Order $\Vert \cdot \Vert_1 $    &Error $\Vert \cdot \Vert_1 $    &Order $\Vert \cdot \Vert_1 $     &Error $\Vert \cdot \Vert_1 $    &Order $\Vert \cdot \Vert_1 $  \\
\hline \hline

 0.1053   &7.94e-03    &       &9.01e-04   &    	&2.09e-04   & \\
 0.0526   &2.08e-03    &1.93   &6.13e-05   &3.88   &4.27e-06   &5.62 \\
 0.0263   &5.22e-04    &1.99   &3.89e-06   &3.98   &7.49e-08 	&5.83\\
 0.0132   &1.29e-04    &2.01   &2.44e-07   &4.00   &1.20e-09 	&5.96\\
 0.0066   &3.08e-05    &2.00   &1.51e-08   &4.00   &1.87e-11 	&6.00\\
 0.0033   &6.16e-06    &2.00   &8.76e-10   &4.00   &2.84e-13 	&6.00\\ 
 \hline
\end{tabular}
\caption{Burgers equation  with initial condition  \eqref{sin_test}, $CFL=0.5$ and $t=0.5$s: $L^1$ errors and accuracy order for LW-CAT2p, $p=1,\, 2,\,3$.}
\label{Tabla_exp07}
\end{center}
\end{table}

\subsection{1D Euler Equation}

We solve the $1D$ Euler equation for gas dynamics 
\begin{align}\label{nonlineareq_euler} 
&\vu_t +  \vf(\vu)_x = 0,
\end{align}
with
\begin{align}
\vu = \left[  \begin{array}{cc}
&\rho \\
&\rho u \\
&E\\
\end{array} \right] ,\quad\vf(\vu) = \left[  \begin{array}{cc}
&\rho u \\
&p + \rho u^2 \\
& u(E+p)\\
\end{array} \right], 
\end{align}
where $\rho$ is the density, $u$ the velocity, $E$  the total energy per unit volume, and $p$  the pressure. 
We assume an ideal gas with the  equation of state, 
\begin{equation}
p(\rho, e) = (\gamma - 1)\rho e,
\end{equation}
being $\gamma$ the ratio of specific heat capacities of the gas taken as 1.4 and $e$ is the internal energy per unit mass given by:
\begin{equation}
E = \rho (e+0.5 u^2).
\end{equation}

We consider the space interval $[0,2]$ with the initial condition:
\begin{align}\label{euler_sin}
\begin{array}{cc}
&\rho(x,0)=0.75+0.5 \sin(\pi x), \\
&\rho u(x,0) =0.25+0.5 \sin(\pi x), \\
&E(x,0)   =0.75+0.5 \sin(\pi x),
\end{array}  
\end{align} and periodic boundary conditions. For this test we take $CFL=0.5$ and $t=0.5$s.  We use  a fine grid with $1400$-point mesh to compute LW-CAT8 as a reference solution.  The results in Table \ref{Tabla_exp_10} support the theoretically obtained accuracy.    

\begin{table}[htbp]
\begin{center}
\begin{tabular}{|l| l l| l l| l l|}
\hline
                 &LW-CAT2         &      &LW-CAT4       &        &LW-CAT6   &    \\
\hline
 $\Delta x$  &Error $\Vert \cdot \Vert_1 $    &Order $\Vert \cdot \Vert_1 $    &Error $\Vert \cdot \Vert_1 $    &Order $\Vert \cdot \Vert_1 $     &Error $\Vert \cdot \Vert_1 $    &Order $\Vert \cdot \Vert_1 $  \\
\hline \hline

0.1053   &3.34e-03    &       &8.57e-04   &    	&5.49e-04   & \\
0.0526   &8.82e-03    &1.92   &9.93e-05   &3.11   &3.53e-05   &4.96 \\
0.0263   &2.28e-04    &1.95   &7.31e-06   &3.76   &1.01e-06 	&5.12\\
0.0132   &5.69e-05    &2.01   &4.81e-07   &3.93   &1.94e-08 	&5.71\\
0.0066   &1.35e-05    &2.07   &3.02e-08   &3.99   &3.21e-10 	&5.92\\
0.0033   &2.71e-06    &2.30   &1.78e-09   &4.08   &4.99e-12 	&6.01\\ 
 \hline
\end{tabular}
\caption{1D Euler equations with initial condition (\ref{euler_sin}), $CFL = 
0.5$ and $t = 0.5$s: $L^1$ errors and accuracy order for LW-CAT2p, $p=1,\, 2,\, 3$.}
\label{Tabla_exp_10}
\end{center}
\end{table}

Finally,  two tests involving  discontinuities are considered: 
\begin{itemize}

\item The Sod Shock tube problem. The initial condition is given by:
$$
(\rho, \rho u, p)
= \left\{
\begin{array}{cl}
\displaystyle (1,0,1) & \mbox {if } x < 0, \\
\displaystyle (0.125,0,0.1) & \mbox {if } x > 0.
\end{array}\right.
$$
 
Here, $x \in [-5,5]$, $CFL=0.5$, $t=1$s, and outflow boundary conditions are considered at both sides. For details of this problem see \cite{SOD19781}. We compare FL-CAT4, WENO5-CAT4, WENO5-LWA5 and WENO5-RK3 using $450$ points. A reference solution is computed  with WENO5-RK3 using a $1 400$-point mesh.  The flux limiter function is computed for every variable.

\begin{figure}[!ht]
	\small
	\setlength{\unitlength}{1mm}
	\centering	
	\begin{picture}(70,65)
	\put(-17,3){\makebox(100,50)[c]{
	\includegraphics[width=16cm,height=8cm]{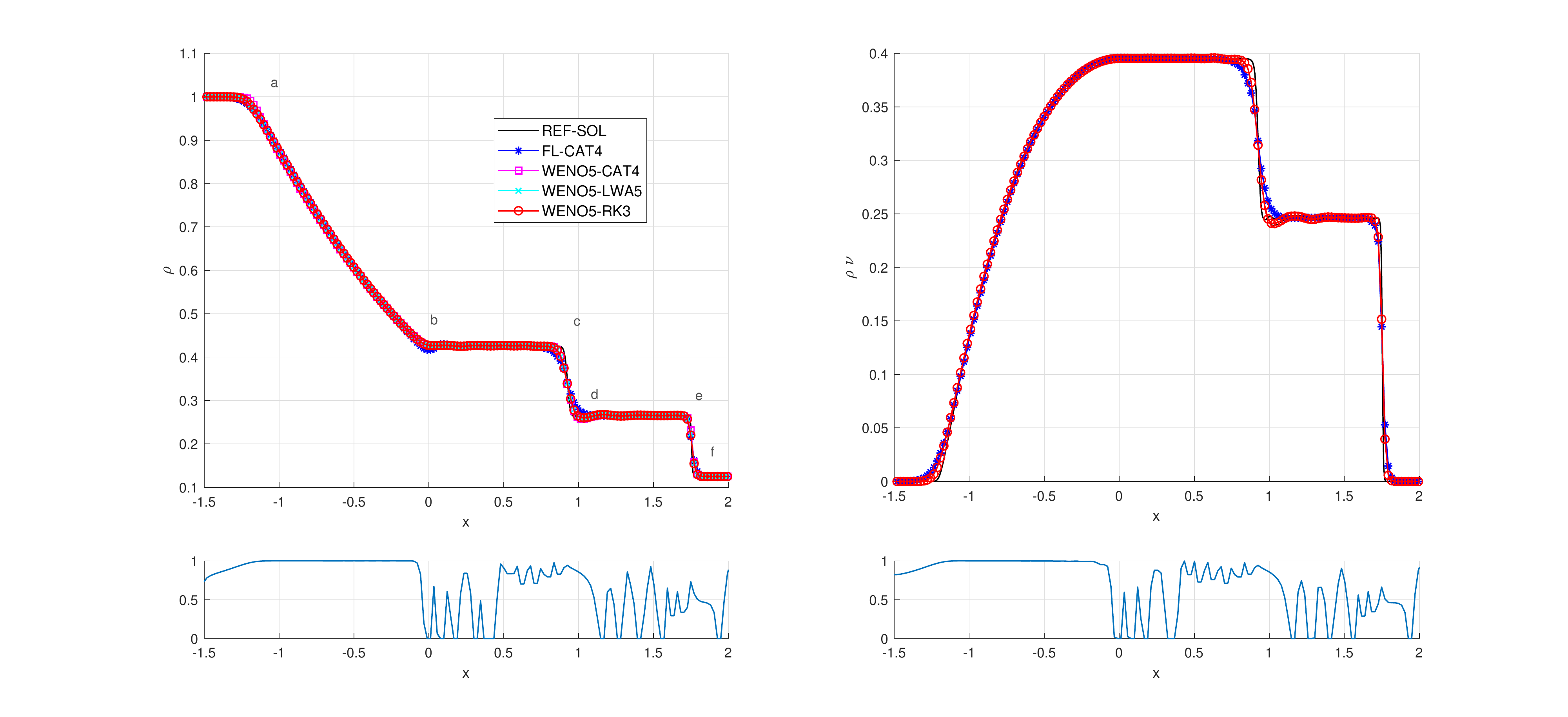}}}
	\end{picture}
	\vspace{1cm}
	\caption{Sod shock tube problem, $CFL=0.5$ and $t=0.5$s.  Left-up: general view of numerical solutions for $\rho$. Left-down:  flux limiter function for $\rho$. Right-up: general view of numerical solutions for $\rho u$. Right-down: flux limiter function for $\rho u$.}
	\label{Num_07}
	\end{figure}

While all numerical solutions show stable and similar values over smooth regions (see Figure \ref{Num_07}), the quality is different in the interest  regions (\textbf{a,b,c,..,f}): an enlarged view of them can be seen in Figure \ref{Num_08}. In the figure we can observe that the solution given by  FL-CAT4 is the most diffusive one (due to the condition  $CFL=0.5$) and that  WENO5-LWA5 and WENO5-RK3 give essentially the same results. WENO5-CAT4 achieves some improvements, specially in  \textbf{a},\textbf{c} and \textbf{d}.

\begin{figure}[!ht]
	\small
	\setlength{\unitlength}{1mm}
	\centering	
	\begin{picture}(70,65)
	\put(-17,3){\makebox(100,40)[c]{
	\includegraphics[width=18cm,height=9cm]{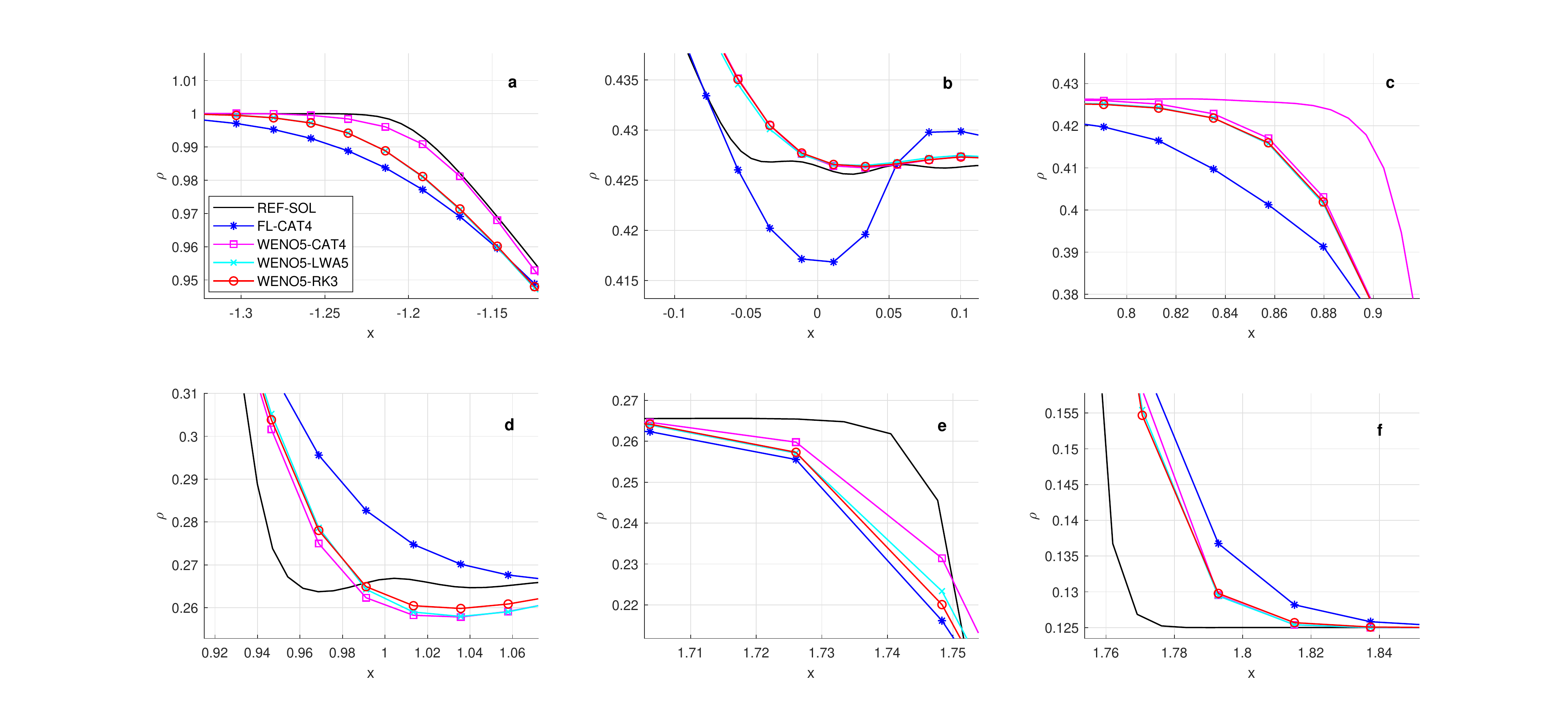}}}
	\end{picture}
	\vspace{1.5cm}
	\caption{Sod shock tube problem, $CFL=0.5$ and $t=0.5$s:  zooms of the numerical results for $\rho$ close to regions  \textbf{a},\textbf{b},\textbf{c},..,\textbf{f}.}
	\label{Num_08}
	\end{figure}

\item The Shu-Osher problem. The initial condition is given by: 
$$
(\rho, \rho u, p)
= \left\{
\begin{array}{cl}
\displaystyle (3.8571,2.6293,10.3333) & \mbox {if } x < -4, \\
\displaystyle (1+0.2 \sin(5x),0,1) & \mbox {if } x > -4.
\end{array}\right.
$$

We consider the space interval $x\in[-5,5]$, $CFL=0.5$ and time $t=1$s. For details see \cite{Shu1989} test 8. We compare FL-CAT4, WENO5-CAT4, WENO5-LWA5 and WENO5-RK3 using $450$- point mesh and a reference solution computed with WENO5-RK3 using a $1 400$-point mesh.  

\begin{figure}[!ht]
	\small
	\setlength{\unitlength}{1mm}
	\centering	
	\begin{picture}(70,65)
	\put(-17,3){\makebox(100,30)[c]{
	\includegraphics[width=18cm,height=9cm]{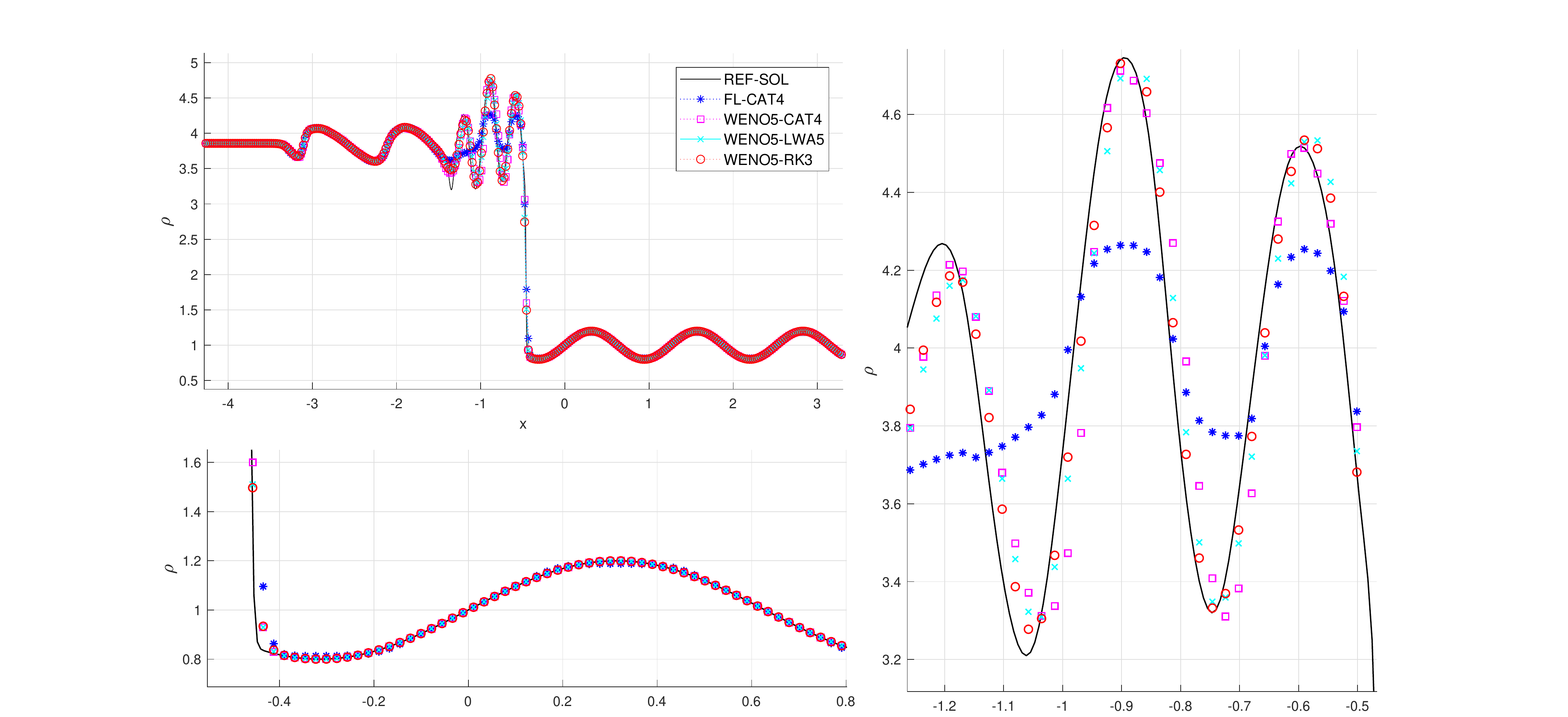}}}
	\end{picture}
	\vspace{3cm}
	\caption{Shu-Osher problem, $CFL=0.5$ and $t=1$s:  Left-up:  general view of numerical solutions for density.  Left-down: enlarged view. Right: enlarged view.}
	\label{Num_09}
	\end{figure}

Again, all solutions are closely similar and near to the reference solution with exception of FL-CAT4.

\end{itemize}

\section{Conclusions}

In this work, first a review of high order Lax-Wendroff methods for the 
linear transport equation has been presented, including the study of 
the order and the $L^2$-stability as well as  the computation and 
properties of the coefficients. Next, an extension to nonlinear 
conservation laws has been introduced with arbitrary even order  $2p$ 
of  accuracy, the so-called Compact Approximate Taylor (CAT) methods. 
Unlike previous applications of Taylor methods to conservation laws, CAT 
methods have $2p + 1$-point centered stencils, like Lax-Wendroff methods 
for linear problems. Moreover, since they inherit the stability 
properties of Lax-Wendroff methods, they are linearly $L^2$-stable under 
a CFL-1 condition.
In order to prevent the spurious oscillations that appear close to 
discontinuities two shock-capturing techniques have been considered: a 
flux-limiter technique (FL-CAT methods) and WENO reconstruction for the 
first time derivative (WENO-CAT methods). We follow  \cite{ZORIO} in 
the second approach.

  These new methods have been compared in a number of test cases with 
WENO-RK methods (Finite Differences WENO reconstructions in space, 
TVD-RK in time) and
  with the WENO-LW methods introduced in  \cite{ZORIO} (Finite 
Differences WENO reconstruction for the first time derivative, 
Approximate Taylor in time).  The linear transport equation, Burgers 
equation,  and the 1D compressible Euler system have been considered.For 
$CFL \leq 0.5$
  all the numerical methods work correctly, and the results obtained 
with all the methods using WENO reconstructions are similar, while the 
FL-CAT method is more diffusive as
  expected. Nevertheless, CAT methods are more expensive in 
computational time and number of operations due to its local character 
(FL-CAT is less expensive than WENO-CAT
  as reconstructions are avoided). However, the extra computational cost 
of CAT methods is compensated by the fact that they still give good 
solutions with CFL values close to 1.

  Future developments include:

  \begin{itemize}

  \item Parallel implementation.

  \item Use of fast WENO reconstructions: see \cite{ZORIO2}  

  \item Order adaptive CAT methods based on smooth indicators.

  \item Application to systems of balance laws.

  \item Extension to multidimensional problems.

\end{itemize}

\section{Acknowledgements}

\begin{figure*}[!ht]
	\small
	\setlength{\unitlength}{1mm}
	\centering	
	\begin{picture}(50,45)
	\put(-55,3){\makebox(100,50)[c]{
\includegraphics[width=4cm,height=3cm]{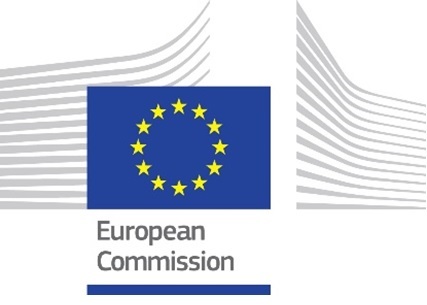}}}
	\end{picture}
	\vspace{-1.2 cm}
	\end{figure*}

“This project has received funding from the European Union’s Horizon 2020 research and innovation program, under the Marie Sklodowska-Curie grant agreement No 642768”.

\end{document}